\documentclass{amsart}
\usepackage{amsmath, bm}
\usepackage{amssymb, framed}
\usepackage{hyperref}
\usepackage{a4wide}
\hypersetup{ colorlinks = true, urlcolor = blue, linkcolor = blue, citecolor = red }
\usepackage{xcolor}
\usepackage{enumitem}
\usepackage{esint}
\makeatletter
\def\namedlabel#1#2{\begingroup
	#2%
	\def\@currentlabel{#2}%
	\phantomsection\label{#1}\endgroup
}
\usepackage{varwidth}
\usepackage{tasks}
\DeclareMathOperator{\loc}{loc}

\newcommand{\RR}{\mathbb{R}}

\newcommand{\pa}{\partial}

\newcommand{ \bb }{ \mathbf{b} }
\newcommand{ \bfu }{ \mathbf{u} }

\newcommand{ \R }{ \mathbb{R} }

\newcommand{ \bA }{ \mathbf{A} }
\newcommand{ \A }{ \mathcal{A} }

\newcommand{ \bV }{ \mathbf{V} }

\newcommand{ \bP }{ \mathbf{P} }

\newcommand{ \bfw }{ \mathbf{w} }
\newcommand{ \bfh }{ \mathbf{h} }
\newcommand{ \zero }{ \mathbf{0} }

\renewcommand{\d}{\mathrm{d}}

\renewcommand{\epsilon}{\varepsilon}
\renewcommand{\phi}{\varphi}
\renewcommand{\le}{\leqslant}
\renewcommand{\ge}{\geqslant}
\renewcommand{\leq}{\leqslant}
\renewcommand{\geq}{\geqslant}
\renewcommand{\div}{\operatorname{div}}

\theoremstyle{plain}
\newtheorem{theorem}{Theorem}[section]
\newtheorem{lemma}[theorem]{Lemma}

\newtheorem{definition}[theorem]{Definition}

\newtheorem{remark}[theorem]{Remark}
\def\Xint#1{\mathchoice
	{\XXint\displaystyle\textstyle{#1}}%
	{\XXint\textstyle\scriptstyle{#1}}%
	{\XXint\scriptstyle\scriptscriptstyle{#1}}%
	{\XXint\scriptstyle\scriptscriptstyle{#1}}%
	\!\int}
\def\XXint#1#2#3{{\setbox0=\hbox{$#1{#2#3}{\int}$}
		\vcenter{\hbox{$#2#3$}}\kern-.5\wd0}}

\def\Yint#1{\mathchoice
	{\YYint\displaystyle\textstyle{#1}}%
	{\YYint\textstyle\scriptstyle{#1}}%
	{\YYint\scriptstyle\scriptscriptstyle{#1}}%
	{\YYint\scriptscriptstyle\scriptscriptstyle{#1}}%
	\!\iint}
\def\YYint#1#2#3{{\setbox0=\hbox{$#1{#2#3}{\iint}$}
		\vcenter{\hbox{$#2#3$}}\kern-.51\wd0}}
\def\longdash{{-}\mkern-3.5mu{-}} 

\def\fiint{\Yint\longdash}

\def\Xint#1{\mathchoice
	{\XXint\displaystyle\textstyle{#1}}%
	{\XXint\textstyle\scriptstyle{#1}}%
	{\XXint\scriptstyle\scriptscriptstyle{#1}}%
	{\XXint\scriptscriptstyle\scriptscriptstyle{#1}}%
	\!\int}
\def\XXint#1#2#3{{\setbox0=\hbox{$#1{#2#3}{\int}$ }
		\vcenter{\hbox{$#2#3$ }}\kern-.6\wd0}}

\def\dashint{\Xint-}

\usepackage{nameref}
\makeatletter
\let\orgdescriptionlabel\descriptionlabel
\renewcommand*{\descriptionlabel}[1]{%
	\let\orglabel\label
	\let\label\@gobble
	\phantomsection
	\edef\@currentlabel{#1}%
	\let\label\orglabel
	\orgdescriptionlabel{#1}%
}
\makeatother
\numberwithin{equation}{section}

\def\Xint#1{\mathchoice
    {\XXint\displaystyle\textstyle{#1}}%
    {\XXint\textstyle\scriptstyle{#1}}%
    {\XXint\scriptstyle\scriptscriptstyle{#1}}%
    {\XXint\scriptscriptstyle\scriptscriptstyle{#1}}%
    \!\int}
\def\XXint#1#2#3{\setbox0=\hbox{$#1{#2#3}{\int}$}
    \vcenter{\hbox{$#2#3$}}\kern-0.5\wd0}
\def\fint{\Xint-}
\def\dashint{\Xint{\raise4pt\hbox to7pt{\hrulefill}}}

\def\XXiint#1#2#3{\setbox0=\hbox{$#1{#2#3}{\iint}$}
    \vcenter{\hbox{$#2#3$}}\kern-0.5\wd0}

\begin{document}
	
\title[Parabolic systems of double phase type]{Partial regularity for parabolic systems \\ of double phase type}

\author{Jihoon Ok}
\address[Jihoon Ok]{Department of Mathematics, Sogang University, 35 Baekbeom-ro, Mapo-gu, Seoul 04107, Republic of Korea.}
\email{jihoonok@sogang.ac.kr}

\author{Giovanni Scilla}
\address[Giovanni Scilla]{Department of Mathematics and Applications ``R. Caccioppoli'', University of Naples Federico II, Via Cintia, Monte S. Angelo, 80126 Naples, Italy}
\email[Giovanni Scilla]{giovanni.scilla@unina.it}

\author{Bianca Stroffolini}
\address[Bianca Stroffolini]{Department of Mathematics and Applications ``R. Caccioppoli'', University of Naples Federico II, Via Cintia, Monte S. Angelo, 80126 Naples, Italy}
\email[Bianca Stroffolini]{bstroffo@unina.it}

\everymath{\displaystyle}

\makeatletter
\@namedef{subjclassname@2020}{\textup{2020} Mathematics Subject Classification}
\makeatother

\begin{abstract}
We study partial regularity for nondegenerate parabolic systems of double phase type, where  the growth function is given by $H(z,s)=s^p+a(z)s^q$, $z=(x,t)\in\Omega_T$, with  $\tfrac{2n}{n+2}<p\le q$ and $a(z)$ a nonnegative $C^{0,\alpha,\frac{\alpha}{2}}$-continuous function for some $\alpha\in(0,1]$. As the main result we prove that if $q< \min \{p+\tfrac{\alpha p }{n+2},  p+1 \}$ the spatial gradient of any weak solution is locally H\"older continuous, except on a set of measure zero. 
\end{abstract}

\keywords{Parabolic double-phase systems, parabolic $\phi$-Laplace systems, partial regularity, gradient estimates}
\subjclass[2020]{35D30, 35K55, 35K65}
\maketitle

\tableofcontents

\section{Introduction}

Double phase energies of the form
\begin{equation}\label{energy}
\int_\Omega H(x,D\bfu)\,\d x\,, \quad H(x,s):= s^p + a(x) s^q\,, \ 1<p<q\,,
\end{equation}
have been originally introduced by Zhikov \cite{Zhikov86,ZKO94} in the context of homogenization theory and Elasticity, with the aim of modeling strongly anisotropic materials, and then studied as an example exhibiting Lavrentiev phenomenon (\cite{Zhikov95,Zhikov97}, see also \cite{espoleomin,FMM04,BDS20}). In this framework, the coefficient $a(\cdot)$ plays a crucial role, since it describes the geometry of a composite made of two different materials exhibiting power-type hardening with exponents $p$ and $q$, respectively. The double phase setting naturally interpolates between two different growth conditions and thus provides a unifying model capturing nonuniform behaviors within a single variational framework. We further refer to \cite{CBZ24,HH21} for the application of the double phase problems to image restoration.

Equations and systems related to \eqref{energy} have been the object of intensive study over the last decade, particularly in connection with regularity theory. For scalar minimizers of \eqref{energy}, sharp and essentially optimal regularity results were obtained in the works of Baroni, Colombo and Mingione \cite{ColomboMingione15,ColomboMingione151,BCM15,BCM18}, where they proved that, under the condition $q\le p+\tfrac{\alpha p}{n}$, the gradient of minimizers is locally Hölder continuous. For further regularity results related to elliptic double phase type problems, we refer to \cite{BaBy25,BBK23,ByunOh17,BOS22,ColomboMingione16,DFDFP,cristianaMin20,cristianaMin23,cristianaOh19,cristianaPalatucci19,OK2017,OKNLA20} and the references therein. We also mention \cite{cristianaMin23b,cristianaMin25,HasOk22,HasLeeOk,HasOk2023,HasOk2022} for regularity results concerning more general energies, including certain double phase type functionals.

Regularity for elliptic systems of the form
\begin{align}\label{system0}
	\operatorname{div}\bA(x, D\bfu) = {\bf 0}\quad \text{in }\ \Omega \,,
\end{align}
where  $\Omega \subset \mathbb{R}^n$ ($n \ge 2$) is a bounded and open set,
$\bfu=(u^1,\dots,u^N): \Omega\to \R^N$ with $N> 1$, 
and the operator $\bA(x, \bm\xi)$ satisfies standard growth and ellipticity conditions, has been actively studied.  It is well known that, without structural assumptions, solutions to systems exhibit only partial regularity, i.e.,   they are regular outside a set of measure zero, see the survey \cite{Mingionedark} and references therein.  An efficient and widely used technique is the so-called $\mathcal{A}$-harmonic approximation, which involves linearizing the functional or partial differential equation around the average gradient in a ball, see \cite{DuGro00, DLSV12}.  This method enables one to locally approximate solutions of nonlinear problems by solutions of linearized equations, facilitating regularity estimates. In contrast to the $\mathcal{A}$-harmonic approximation, which is suitable when the gradient is bounded away from zero, the regime where the gradient is small requires approximation schemes that reflect the degeneracy or singularity of the underlying operator, see \cite{DuzMin041, DieStrVer12, CelOk20}, respectively for the power case and the Orlicz one.

More recently, partial regularity for nondegenerate systems of the type \eqref{system0} 
with operators $\bA(x, \bm\xi)$ associated with the double phase function $H(x,s)=s^p+a(x)s^q$ has been investigated in \cite{OKJFA18,OKNLA18,SciStroffofluids}. 
Here “nondegenerate’’ refers to the condition $|D_{\bm \xi}\bA(x,0)|\sim 1$.  Those works deal only with the superquadratic case, namely $p\ge 2$.
The general (possibly degenerate) case of double phase systems 
was recently treated in \cite{OkScStdbph}, 
where partial Hölder regularity of the gradient of weak solutions was obtained without the superquadratic assumption, 
thus providing a unified approach independent of the exponent $p$.

For the corresponding parabolic systems
\begin{align}\label{system1}
	\bfu_t - \operatorname{div}\bA(z, D\bfu) = {\bf 0}
	\quad \text{in } \Omega_T,
\end{align}
where $z=(x,t)$ and $\Omega_T=\Omega\times(0,T)$ is a space-time cylinder in $\R^{n+1}$,  partial regularity in the $p$-growth case has been studied in \cite{DuzMinSte11,BoDuMin13}. Moreover, in recent years, partial regularity in the Orlicz growth setting has also been investigated in \cite{FILV23,OkSciStroffo2024}.

In contrast to the elliptic setting, where the picture of regularity results is rich and appears to be almost completely understood, the corresponding parabolic theory for double phase equations and systems remains much less developed. 
 Chlebicka, Gwiazda and Zatorska-Goldstein \cite{ChGwZa19} proved the existence of weak solutions to parabolic equations with Musielak--Orlicz growth, and Kim, Kinnunen and S\"arki\"o \cite{KKSA} proved existence and uniqueness results for weak solutions to parabolic systems with superquadratic double phase growth. 
Concerning regularity theory, Kim, Kinnunen, Moring and S\"arki\"o \cite{K1,K2,KKM23,KSA24} have established higher integrability for the gradient as well as Calder\'on--Zygmund type estimates for weak solutions. These results were obtained under the conditions
\begin{equation}\label{eq:pqKinnunen}
q \le p+\frac{2\alpha}{n+2} \quad \text{when } p \ge 2, 
\qquad 
q \le p+\frac{\alpha}{n+2}\,\frac{p(n+2)-2n}{2} 
\quad \text{when } \frac{2n}{n+2} < p < 2 \,.
\end{equation}
Moreover, in \cite{KKM23} an additional technical assumption $D\bfu \in L^q(\Omega_T)$ was required; however, this can be removed by applying a double phase parabolic Lipschitz truncation, see \cite{KKSA}. We further refer to \cite{KMSA} for results on the H\"older continuity of bounded weak solutions to parabolic equations, 
\cite{CW24} for parabolic double phase obstacle problems, and  \cite{BoDuMarc13,Singer} for related work on parabolic $(p,q)$-growth problems.

To the best of our knowledge, partial regularity for parabolic systems of double phase type has not yet been explored. In this paper, we address this problem for homogeneous nondegenerate systems of the form
\eqref{system1}, where the nonlinearity ${\bf A}:\Omega_T \times \mathbb{R}^{N\times n}\to \mathbb{R}^{N\times n}$ satisfies the following growth and ellipticity conditions with respect to the function $H(z,s)$: 
\begin{equation}
|{\bf A}(z,\bm\xi)|+ (1+|\bm\xi|)\,|D_{\bm\xi} {\bf A}(z,\bm\xi)| \leq L\, H'(z, 1+|\bm\xi|)\,, 
\tag{A1}
\label{eq:1.8ok1}
\end{equation}
\begin{equation}
[D_{\bm\xi} {\bf A}(z,\bm\xi)\bm\lambda]: \bm\lambda \geq \nu\, H''(z, 1+|\bm\xi|)\,|\bm\lambda|^2\,,
\tag{A2}
\label{eq:1.8ok2}
\end{equation}
for every $z=(x,t)\in\Omega_T$, $\bm\xi\in\mathbb{R}^{N\times n}$, and $\bm\lambda \in \mathbb{R}^{N\times n}$, and for some $0<\nu\leq L$.
Here, ``$:$'' denotes the Euclidean inner product in $\mathbb{R}^{N\times n}$, and  the function $H:\Omega_T \times [0,\infty)\to [0,\infty)$ is the so-called double phase function defined by
\begin{equation}
H(z,s):= s^p + a(z)\,s^q\,, \quad \text{where } \ \frac{2n}{n+2} < p \le q, \ \ 0 \le a(z) \le L \ \ \text{for } z\in \Omega_T\,.
\label{eq:H}
\end{equation}
We further assume that $a\in C^{0,\alpha,\frac{\alpha}{2}}(\Omega_T)$ for some $\alpha \in (0,1]$, that is, 
\begin{equation}
|a(z_1)-a(z_2)| \leq L \big(|x_1-x_2|^\alpha + |t_1-t_2|^{\frac{\alpha}{2}}\big)
\label{eq:holdercond}
\end{equation}
for every $z_1=(x_1,t_1)$ and $z_2=(x_2,t_2)$ in $\Omega_T$, and  
\begin{equation}
q< \min \Big\{p+\frac{\alpha p }{n+2}, \, p+1 \Big\}\, .
\label{eq:pq}
\end{equation}
Moreover,  $H'(z,s)$ and $H''(z,s)$ denote the first and second derivatives with respect to the $s$-variable, respectively.
Note that in the region where $a(z)=0$, we have $H(z,s)=s^p$, so that $H$ corresponds to the $p$-phase, while in the region where $a(z)>0$, the function $H$ corresponds to the $(p,q)$-phase.  
We also note that \eqref{eq:1.8ok2} implies the following monotonicity property: for every $z\in\Omega_T$ and $\bm\xi_1, \bm\xi_2\in \mathbb{R}^{N\times n}$,
\begin{equation}
\big({\bf A}(z,\bm\xi_1) - {\bf A}(z,\bm\xi_2)\big) \, : \, (\bm\xi_1-\bm\xi_2) 
\geq \tilde{\nu}\, H''(z, 1+ |\bm\xi_1|+|\bm\xi_2|)\,|\bm\xi_1-\bm\xi_2|^2
\label{eq:1.9ok}
\end{equation}
for some $\tilde \nu>0$ depending only on $p$, $q$, and $\nu$.

For system \eqref{system1}, we establish the partial $C^{1,\alpha}$-regularity theory; see Section~\ref{sec:mainresult} for the main theorem. We outline the novelty of our result and the main ideas of the proof. Our analysis is focused on non-degenerate systems, which represent the initial step in the study of partial regularity theory. However, even in this context, the proof of partial regularity necessitates the introduction of significantly novel concepts.
In the nondegenerate setting, a natural condition on the exponents $p$ and $q$ in the definition of $H$ might be
\begin{equation}\label{eq:pq=}
q \le p+\frac{p\alpha}{n+2}\,,
\end{equation}
as this allows one to avoid the use of intrinsic parabolic cylinders. We point out that, even in the elliptic case, in order to obtain regularity in the borderline situation $q=p+\tfrac{p\alpha}{n+2}$, higher integrability plays an essential role. In the parabolic setting, it seems that the conditions in \eqref{eq:pqKinnunen} are essentially required to obtain higher integrability, and are more restrictive than \eqref{eq:pq=}. Therefore, assuming the strict inequality in \eqref{eq:pq=} is quite natural. 

A Caccioppoli type inequality was obtained in \cite{KKSA} for weak solutions $\bfu$ to the double phase parabolic systems satisfying the natural integrability condition $H(\cdot,|D\bfu|)\in L^1_{\loc}(\Omega_T)$, assuming the strong condition \eqref{eq:pqKinnunen}. On the other hand, in this paper,
we are able to prove it under the optimal condition \eqref{eq:pq=} by means of a mollification argument of convolution type in both space variable $z$ and the time variable $t$, originally devised for the Euclidean metric in \cite{HH} and later adapted to the parabolic setting in \cite{HasOk25}.

Thanks to our recent $\mathcal{A}$-caloric approximation result in \cite{OkSciStroffo2024}, we can define the \emph{excess functional} directly in terms of the gradient $D\bfu$ of a weak solution $\bfu$ to \eqref{system1}:
\begin{equation}\label{eq:excessintro}
\Phi(z_0,\rho,\bfu) = \fint_{Q_{\rho}(z_0)} H^-_{1+|(D\bfu)_{z_0,\rho}|}
   \left(\big|D\bfu -(D\bfu)_{z_0,\rho}\big|\right)\,\d z + \rho^{\frac{\beta_1}{2}}\,,
\end{equation}
where $H^-_{\sigma}$ denotes the shifted $N$-function of $H^-(s):= s^p+\inf_{z\in Q_\rho(z_0)}a(z) s^q$ and $\beta_1$ is defined in \eqref{def:beta1}. 
Then a linearization procedure combined with $\mathcal{A}$-caloric approximation yields local comparison estimates for $D\bfu$ with the gradient of a suitably chosen smooth $\mathcal{A}$-caloric function. 
We note that the additional assumption \( q < p+1 \) is imposed  in order to apply the \(\mathcal{A}\)-caloric approximation to the function \( H^- \), rather than \( H \). This condition could, in principle, be relaxed if a higher integrability result for functions of the form \( H_{1+|{\bf Q}|}(z,|D\bfu - {\bf Q}|) \) were to be established. However, as previously noted, this would likely require the conditions in \eqref{eq:pqKinnunen}. On the other hand, the condition \( q < p+1 \) is less restrictive than \eqref{eq:pqKinnunen}, and naturally arises in the regularity theory for parabolic problems with \((p,q)\)-growth; see, for instance, \cite{BoDuMarc13,BoDuMarc13b}.

We then establish a decay estimate for the excess functional \eqref{eq:excessintro} in Lemma~\ref{Lem:nondegenerate1}. Finally, by iterating this decay at smaller scales, as shown in Lemma~\ref{Lem:iteration}, we obtain the partial regularity result in the standard manner, as stated in Section~\ref{sec:mainthm}: the weak solution is locally Hölder continuous, except on a set of measure zero; see Theorem~\ref{thm:main-thm}.


As a final remark, we highlight an important aspect of our approach. In most papers on double phase problems, different methods have been adopted for the $p$-phase and the $(p,q)$-phase. However, in this paper, the arguments of approximation, comparison, and iteration do not distinguish between the phases. This allows us to provide simpler and more unified proofs.

\emph{Outline of the paper:} The paper is organized as follows.  In Section~\ref{sec:mainresult}, we present our main result, Theorem~\ref{thm:main-thm}. In Section~\ref{sec:preliminaries}, we fix the necessary notation, and recall some basic facts and elementary inequalities about Orlicz and the generalized Orlicz function $H$. In Section~\ref{sec:caccioppoli}, we obtain a Caccioppoli type estimate by using a new mollification argument. 
In Section~\ref{sec:approxandlin}, we collect two main ingredients in order to get the decay estimates: we recall regularity estimates for the solutions of systems with constant coefficients and the $\mathcal{A}$-caloric approximation, see Section~\ref{subsec:CZ}, and then we setup the linearization, Section~\ref{sec:linearization}. Section~\ref{sec:decayestimate} provides decay estimates for the excess functional (Lemma \ref{Lem:nondegenerate1}), which can be iterated on shrinking parabolic cylinders (Lemma \ref{Lem:iteration}).  Finally, in Section~\ref{sec:mainthm}, we prove our main result.

\subsection{Main result}\label{sec:mainresult}

We introduce the main result of the paper. First, we list the required additional assumptions on ${\bf A}(z,\bm\xi)$: we assume the following continuity condition on $z\mapsto {\bf A}(z,\bm\xi)$  
\begin{equation}\label{eq:Acontinuity}
\begin{split}
|{\bf A}(z_1,\bm\xi) - {\bf A}(z_2,\bm\xi)| & \leq L (|x_1-x_2| + \sqrt{|t_1-t_2|})^{\beta_0} \left(H^\prime (z_1, 1+|\bm\xi|) + H^\prime (z_2, 1+|\bm\xi|)\right) \\
& \,\,\,\,\,\, + L |a(z_1)-a(z_2)| (1+|\bm\xi|)^{q-1}\,
\end{split}
\tag{A3}
\end{equation} 
for all $z_1=(x_1,t_2)$ and $z_2=(x_2,t_2)$ in $\Omega_T$ and for some $\beta_0\in(0,1)$,  
and the off diagonal condition on $z\mapsto D_{\bm\xi}{\bf A}(z,\bm\xi)$  
\begin{equation}\label{offdiagonal}
|D_{\bm\xi}\bA(z,\bm\xi_1) - D_{\bm\xi}\bA(z,\bm\xi_2)| \le L \left(\frac{|\bm\xi_1-\bm\xi_2|}{1+|\bm\xi_1|}\right)^{\gamma}H''(z, 1+|\bm\xi_1|)
\ \text{ for some }\ \gamma\in(0,1)\,,
\tag{A4}
\end{equation}
for  all $\bm\xi_1,\bm\xi_2\in \R^{N\times n}$ with $2|\bm\xi_1-\bm\xi_2|\le 1+|\bm\xi_1|$.

We then specify the notion of weak solution to \eqref{system1}. 
\begin{definition}
A function $\bfu=(u^1,u^2,\dots,u^N) \in C_{\loc}(0,T; L^2_{\loc}(\Omega,\mathbb{R}^N)) \cap L^p_{\loc}(0,T; W^{1,p}_{\loc}(\Omega,\mathbb{R}^N))$ with $H(\cdot,|D{\bf u}|) \in L^1_{\loc}(\Omega_T)$ is said to be a  \emph{weak solution} to the system \eqref{system1} if it satisfies the following weak form of \eqref{system1}: 
\begin{equation}
-\iint_{\Omega_T} \bfu \cdot {\bm\zeta}_t \,\d z + \iint_{\Omega_T}   \bA(z , D\bfu) : D{\bm\zeta} \, \d z =0 
\quad  \text{for all }\ {\bm\zeta}\in C^\infty_{\mathrm c}(\Omega_T,\R^N)\,,
\label{eq:weakformul}
\end{equation}
where  ``$\cdot$" is the Euclidean inner product in $\R^N$. 
\end{definition}

Finally, we are in position to state the main result of the paper. 

\begin{theorem}\label{thm:main-thm}
Let $H:\Omega_T\times[0,\infty)\to[0,\infty)$ be defined as in \eqref{eq:H} complying with \eqref{eq:holdercond} and \eqref{eq:pq}, and  ${\bf A}:\Omega_T \times\R^{N\times n}\to \R^{N\times n}$ comply with \eqref{eq:1.8ok1}--\eqref{offdiagonal}.
If $\bfu \in C_{\loc}(0,T; L^2_{\loc}(\Omega,\mathbb{R}^N)) \cap L^p_{\loc}(0,T; W^{1,p}_{\loc}(\Omega,\mathbb{R}^N))$ with $H(\cdot,|D{\bf u}|) \in L^1_{\loc}(\Omega_T)$ is a weak solution to \eqref{system1}, then there exist $\beta=\beta(n,p,q,\alpha,\beta_0)\in(0,1)$ and an open subset $\Omega_0 \subset \Omega_T$ such that $|{\Omega_T \setminus \Omega_0} |  = 0$ and
\begin{equation*}
 D\bfu \in C^{0, \beta, \frac{\beta}{2}}_{\rm loc}\left(\Omega_0 ;\R^{N\times n}\right)  \,.
\end{equation*}
Moreover, $\Omega_T \setminus \Omega_0\subset \Sigma_1\cup\Sigma_2$ where
\begin{equation*}
\begin{split}
&\Sigma_1:=\left\{z_0\in\Omega_T:\,\, \liminf_{r \to 0^+} \fiint_{Q_r(z_0)} |{\bf V}_{H^-_{Q_r(z_0)}}(D\bfu) - ({\bf V}_{H^-_{Q_r(z_0)}}(D\bfu))_{z_0,r}  |^2 \, \mathrm{d}z>0\right\}\,,\\
&\Sigma_2:=\left\{z_0\in\Omega_T:\,\, \limsup_{r \to 0^+}  | (D\bfu)_{z_0,r} |  = \infty\right\}\,.
\end{split}
\end{equation*} 
where $H^-_{Q_r(z_0)}(s)$ and  ${\bf V}_{H^-_{Q_r(z_0)}}({\bf P})$ are defined as in \eqref{eq:Hpmdef}  and \eqref{Vfunction} with $\varphi(s)=H^-_{Q_r(z_0)}(s)$, respectively.
\end{theorem}

\begin{remark}
Notice that the set $\Sigma_1$ has measure zero, since it is contained in the set of (parabolic) non-Lebesgue points of ${\bf V}_{H}(\cdot,D\bfu)$ (see \eqref{VfunctionH} for the definition of ${\bf V}_H$) in $L^2$-space. For details, see \cite[Remark 1]{OkScStdbph}.
\end{remark}

\section{Preliminaries}\label{sec:preliminaries}

\subsection{Basic notation} 
A space-time cylinder in $\RR^{n+1}$ is denoted by 
\[
Q_{r,\tau}(z_0)=B_{r}(x_0)\times (t_0-\tau,t_0+\tau),
\quad \text{where }\ z_0=(x_0,t_0)\in \R^{n+1} \ \text{ and }\ r,\tau>0 \,.
\]
In particular, if $\tau=r^2$, we write $Q_r(z_0):= Q_{r,r^2}(z_0)$ and call it a parabolic cylinder.

Let  ${\bm \ell}: \R^n \to \R^N$ be 
any
linear map of the form
\begin{equation}\label{def_l}
{\bm \ell}(x) : = \bP (x-x_0) + \bb ,
\quad x\in \R^n,
\end{equation}
where $\bP \in \R^{N\times n}$, $x_0\in \R^n$ and $\bb \in \R^N$. If $\bfu$ is a weak solution to \eqref{system1}, we  set
\begin{equation}\label{def_ul}
\bfu_{\bm \ell}(z):= \bfu(z)-{\bm \ell}(x)\,, \quad z=(x,t)\in \Omega_T\,.
\end{equation}
Let $H:\Omega_T\times [0,\infty)\to [0,\infty)$ be  given in \eqref{eq:H} with \eqref{eq:holdercond}. For $Q_r(z_0)\Subset \Omega_T$,  let $z^-_{z_0,r}, z^+_{z_0,r}\in \overline{Q_r(z_0)}$ be such that
\begin{equation}\label{zpm}
a^-_{z_0,r}:=a(z^-_{z_0,r}) = \inf_{z\in Q_r(z_0)}a(z)
\quad\text{and}\quad
a^+_{z_0,r}:=a(z^+_{z_0,r}) = \sup_{z\in Q_r(z_0)}a(z)\,. 
\end{equation}
Then we write
\begin{equation}
\label{eq:Hpmdef}
H^-_{Q_r(z_0)}(s):=H(z^-_{z_0,r},s)
\quad\text{and}\quad
H^+_{Q_r(z_0)}(s):=H(z^+_{z_0,r},s) \,. 
\end{equation}

\subsection{Orlicz functions and operators}

We recall basic notation and  properties about Orlicz functions. The following definitions and results can be found, e.g., in \cite{Kras, Kufn}.

In this paper, $\phi:[0,\infty)\to[0,\infty)$ is always an $N$-function; that is, $\phi(0)=0$, there exists a right continuous derivative $\phi'$ of  $\phi$ such that $\phi'(0)=0$, $\phi'$ is increasing, and $\phi'(s)>0$ when $s>0$. 
We say that $\phi$ satisfies the $\Delta_2$ condition, denoted by $\Delta(\phi)<\infty$, if there exists a positive constant $K =: \Delta(\phi)$ such that $\phi(2s) \le K\phi(s)$ for all $s>0$.   
%
The conjugate function of $\phi$ is defined as
$$
\phi^*(\tilde s):= \sup_{\tilde s\geq 0}\, (s\tilde s - \phi(s))\,.
$$
From the definition of $\phi^*$, the following Young's inequality  
 $$
 s\tilde s  \leq \phi(\tilde s) +  \phi^\ast(s)\,,\quad s,\tilde s\geq0\,,
$$
holds true. From now on, we always assume that $\phi$ and $\phi^*$ satisfy the $\Delta_2$ condition and  this is 
indicated by $\Delta(\phi,\phi^*)<\infty$, where $\Delta(\phi,\phi^*)$ denotes  
the constants $\Delta(\phi)$ and $\Delta(\phi^*)$.
We note that the exact value of $\phi^*$ is not always explicitly computable and instead the estimate
\begin{equation}
\phi^*\left(\frac{\phi(s)}{s}\right)\sim \phi^*(\phi'(s)) \sim \phi(s)
\label{eq:hok2.4}
\end{equation}
will often be useful in computations (see \cite[Theorem~2.4.10]{HH}). Moreover, if $\phi\in C^1([0,\infty))$ and $\phi'$ is strictly increasing, then \eqref{eq:hok2.4} implies that  
 \begin{equation}
(\phi^*)' (s) \sim (\phi')^{-1}(s), \quad s \ge 0.
\label{eq:hok2.41}
\end{equation}

We recall the notion of ``\textit{almost convexity}" and a Jensen type inequality (see \cite[Lemma 4.3.1]{HH}).
\begin{lemma}\label{lem:Jensen}
If $\psi:[0,\infty)\to [0,\infty]$ is increasing with $\psi(0)=0$ and $\tfrac{\psi(s)}{s}$ is almost increasing; i.e., $\frac{\psi(s)}{s}\le L\frac{\psi(\tilde s)}{\tilde s}$ for every $0< s \le \tilde s$ with constant $L\geq 1$, then the following Jensen's type inequality holds:
\[
\psi\bigg( \frac{1}{L^2} \fint_U |f|\, dz\bigg) \le \fint_U \psi(|f|)\, dz.
\]
\end{lemma}

Now we further assume for $\phi$ that 
$\phi\in C^1([0,\infty))\cap C^2((0,\infty))$ and
\begin{equation}\label{pq}
1 < p \le \frac{s\phi''(s)}{\phi'(s)} +1 \le q \quad \text{for all }\ s>0.
\end{equation}
Note that this implies
$$
1< p \le  \frac{s \phi'(s)}{\phi(s)} \le q \quad \text{for all }\ s>0\,,
$$
and hence the $\Delta_2$ conditions of $\phi$ and $\phi^*$.


For an $N$-function  $\phi$ and for $\sigma\ge0 $,  we define the \textit{shifted $N$-function} $\phi_\sigma$ of $\phi$ by
$$
\phi_{\sigma}(s):=\int_{0}^{s} \frac{\phi'(\sigma+\tilde{s}) \tilde{s}}{\sigma+\tilde{s}} \, \d \tilde{s}\,,  \quad \text{that is, }\ \phi_\sigma'(s) =\frac{\phi'(\sigma+s)}{\sigma+s}s\,.
$$
%
The condition \eqref{pq} implies that
\begin{equation}
\label{phiapq}
\tilde p:=\min\{p,2\} \le \frac{s \phi_\sigma''(s)}{\phi'_\sigma(s)} +1 \le \max\{q,2\}=:\tilde q
\quad \text{for all }\ s>0 \ \text{ and }\ \sigma\ge0\,. 
\end{equation}
Moreover, we have the following relations (see, e.g., \cite[Proposition 2.3]{CelOk20} and \cite{DieEtt08,DieStrVer09}), which hold uniformly with respect to $\sigma\geq0$:
\begin{align}
&\phi_\sigma(s) \sim \phi'_\sigma(s)\,s\,; \label{(2.6a)} \\
&\phi_\sigma(s) \sim \phi''(\sigma+s)s^2\sim\frac{\varphi(\sigma+s)}{(\sigma+s)^2}s^2\sim \frac{\varphi'(\sigma+s)}{\sigma+s}s^2\,;\label{(2.6b)}\\
& \phi(s) \le \phi(\sigma+s)\sim \phi_\sigma(s)+\phi(\sigma)\,.\label{eq:approx}
\end{align}

We define vector valued functions $\bV_{\phi}:\R^{N\times n}\to \R^{N\times n}$
by
\begin{equation}\label{Vfunction}
\bV_\varphi({\bf Q}):= \sqrt{\frac{\phi_1'(|{\bf Q}|)}{|{\bf Q}|}} {\bf Q} =\sqrt{\frac{\phi'(1+|{\bf Q}|)}{1+|{\bf Q}|}} {\bf Q}.
\end{equation}
In particular, we write
$$
\bV_p({\bf Q}):= \bV_\varphi({\bf Q})
\quad \text{when } \ \phi(s)=s^p\,.
$$
Then we recall equivalent relations in \cite[Lemmas 3 and  20]{DieEtt08} and \cite[Lemma~3.1]{DieSchSch19} with $\phi_1$ in place of $\phi$:
\begin{equation}\label{monotonicity1}
\frac{\phi'(1+|{\bf P}| + |{\bf Q}|)}{1+|{\bf P}| + |{\bf Q}|}  |{\bf P}-{\bf Q}|^2 \sim |\bV_{\phi}( {\bf P})-\bV_{\phi}( {\bf Q})|^2   \sim \varphi_{1+|{\bf Q}|}(| {\bf P} -  {\bf Q}|) \,,
\end{equation}
\begin{equation}\label{DE08lem20}
\frac{\phi'(1+|{\bf P}| + |{\bf Q}|)}{1+|{\bf P}| + |{\bf Q}|} \sim \int_0^1 \frac{\phi'(1+|\tau {\bf P} + (1-\tau) {\bf Q}|)}{1+|\tau {\bf P} + (1-\tau) {\bf Q}|} \, \d \tau \,.
\end{equation}
Moreover, by the same proof of \cite[Lemma A.2]{DieKaSch12}, we have that for every ${\bf g}\in L^\phi(Q_r; \R^{N\times n})$,
\begin{equation}\label{Vintequivalent}
\fint_{Q_r} |\bV({\bf g}) - (\bV({\bf g}))_{Q_r}|^2\,\d z \sim   \fint_{Q_r} |\bV({\bf g}) - \bV(({\bf g})_{Q_r})|^2\,\d z \,.
\end{equation}
 {Note that all constants concerned with the relation $\sim$ and $c$  in above depend only on $p$ and $q$.}



\subsection{Generalized Orlicz function $H$}

For basic definitions and properties about generalized Orlicz functions and the associated spaces, we refer to the monograph \cite{HH}.

Let $H:\Omega_T\times [0,\infty)\to [0,\infty)$ be the function given in \eqref{eq:H}. Note that for each $z\in \Omega_T$, $\phi(s):= H(z,s)$  is an $N$-function satisfying that $\Delta_2(\varphi,\varphi^*)<\infty$.  With this $H$, we define the generalized Orlicz space $L^H(\Omega_T,\R^N)$ as the set of all measurable functions $f:\Omega_T\to\R^N$ such that
$$
\iint_{\Omega_T} H(z,|f(z)|)\, \d z= \int_0^T\int_\Omega H(x,t,|f(x,t)|)\, \d x \, \d t < \infty, 
$$ 
endowed with the usual Luxembourg type norm
\[
\|f\|_{L^H(\Omega_T)}:= \inf\left\{\lambda >0 : \iint_{\Omega_T} H\left(z,\frac{|f(z)|}{\lambda}\right)\, \d z \le 1 \right\}\,.
\]
%
We define
\begin{equation}\label{VfunctionH}
\bV_H(z,{\bf Q}):= \sqrt{\frac{H'(z,1+|{\bf Q}|)}{1+|{\bf Q}|}} {\bf Q}\,, \quad z\in\Omega_T\,.
\end{equation}
Then it follows from \eqref{monotonicity1} that for every  $z\in\Omega_T$,
\begin{equation}\label{monotonicity1H}
\begin{aligned}
|\bV_H(z,{\bf P})-\bV_H(z, {\bf Q})|^2  & \sim H_{1+|{\bf Q}|}(z,| {\bf P} -  {\bf Q}|) \\
& \sim \phi_{p,1+|{\bf Q}|}(| {\bf P} -  {\bf Q}|) + a(z) \phi_{q,1+|{\bf Q}|}(| {\bf P} -  {\bf Q}|) \\
& \sim |\bV_p({\bf P})-\bV_p( {\bf Q})|^2 + a(z) |\bV_q({\bf P})-\bV_q( {\bf Q})|^2\\
& \ge |\bV_p({\bf P})-\bV_p( {\bf Q})|^2 \,,
\end{aligned}\end{equation}
where $\phi_{p,\sigma}(s)$ is the shifted $N$ function of the power function $s^p$ with shift $\sigma\ge 0$.
Moreover, combining \eqref{DE08lem20} and \eqref{eq:1.8ok1}, following \cite[Lemma~3.1]{DieSchSch19} we can prove that for every $z\in\Omega_T$,
\begin{equation}
|{\bf A}( z, {\bf P})-{\bf A}(z, {\bf Q})|  \sim H'_{1+|{\bf Q}|}(z,|{\bf P} - {\bf Q}|)\,.
\label{eq:(3.4)}
\end{equation}
We also notice that for the shifted $N$-functions $H^\pm_\sigma$ and the corresponding $V_{H^\pm_\sigma}$ function the equivalence \eqref{Vintequivalent} still holds true. This is  related to the set of Lebesgue points.

We end the section with the following parabolic 
Sobolev-Poincar\'e 
type inequality for shifted $N$-functions of $H$.

\begin{lemma}[Sobolev-Poincar\'e inequality]\label{thm:sob-poincare}
Let $H$ be defined as in \eqref{eq:H} with \eqref{eq:holdercond} and \eqref{eq:pq=},
and let $\sigma\ge 0$ and  $0<r\leq1$. 
Then there exists $\theta=\theta(n,p,q)\in(0,1)$ such that for $Q_r=Q_r(z_0)\subset\Omega_T$ with $z_0=(x_0,t_0)$, if ${\bf w}\in L^1(t_0-r^2,t_0+r^2;W^{1,1}(B_r(x_0);\R^N))$ with $\| H^-_\sigma(|D{\bf w}|)\|_{L^1(Q_r)}<\infty $
 satisfies that   $\bfw(\cdot,t)= {\bf 0}$ on $\partial B_r(x_0)$ in the trace sense for a.e. $t\in (t_0-r^2,t_0+r^2)$, then we have 
$$
\fiint_{Q_r} H^+_{\sigma}\left( \frac{|{\bf w}|}{r}\right)\,\mathrm{d}z\leq c_P
\left(\fiint_{Q_r} H^-_\sigma(|D{\bf w}|)^{\theta}\,\mathrm{d}z\right)^\frac{1}{\theta} + c_P \left(r^{\alpha-\frac{(q-p)(n+2)}{p}}+ r^\alpha \sigma^{q-p}\right)\sigma^p\,,
$$
for some $c_P=c_P(n,p,q,L)\geq 1$, where $H^\pm_\sigma$ denote the shifted $N$-functions of $H^\pm_{Q_r}$ (see \eqref{eq:Hpmdef}) with shift $\sigma\ge 0$.
 \end{lemma}

\begin{proof}
Recall $z^\pm_{0,r}$ from \eqref{zpm} and,  for ease of notation, set $z^\pm:= z^\pm_{0,r}$. First, by \eqref{(2.6b)} and \eqref{eq:holdercond} we get 
$$\begin{aligned}
H_{\sigma}^+(s) & \sim  (\sigma+s)^{p-2}s^2 + a(z^+)(\sigma+s)^{q-2}s^2 \\
& \lesssim (\sigma+s)^{p-2}s^2 + a(z^-)(\sigma+s)^{q-2}s^2 + r^\alpha(\sigma+s)^{q-2}s^2\\
& \lesssim  H^-_{\sigma}(s) + r^\alpha(\sigma^q+ s^q)\,.
\end{aligned}$$
Now, for each fixed time $t\in (t_0-r^2,t_0+r^2)$, we can estimate, by using Sobolev-Poincar\'e inequalities for the shifted $N$-function $H^-_{\sigma}(s)$ (see, e.g., \cite[Theorem 7]{DieEtt08}) and for the function $s\mapsto s^q$, 
$$\begin{aligned}
\fiint_{Q_r} H_{\sigma}\left(z, \frac{|{\bf w}|}{r}\right)\,\mathrm{d}z
&\lesssim \fiint_{Q_r} H^-_{\sigma}\left(\frac{|{\bf w}|}{r}\right)\,\mathrm{d}z + r^\alpha \fiint_{Q_r} \left(\frac{|{\bf w}|}{r}\right)^q\,\mathrm{d}z+r^\alpha \delta^q \\
&\lesssim  \left(\fiint_{Q_r} H^-_{\sigma}(|D{\bf w}|)^{\tilde\theta}\,\mathrm{d}z\right)^\frac{1}{\tilde\theta}+ r^\alpha\left(\fiint_{Q_r} |D{\bf w}|^{q_*}\, \mathrm{d}z\right)^{\frac{q}{q_*}} + r^\alpha \sigma^q,
\end{aligned}$$
where $\tilde \theta\in(0,1)$ depends on $n,p,q$ and $q_*:=\max\{1, \tfrac{nq}{n+q}\} \le \max\{1, q\tfrac{n}{n+1}\}<p$. Set  $\theta:=\max\{\tilde \theta,q_*/p\}\in(0,1)$. 
Note that, by \eqref{eq:approx} applied to the $N$-function $s\mapsto s^p$, we have
\begin{equation}
s^p \lesssim (\delta+\sigma)^{p-2}s^2 + \sigma^p \lesssim H^-_\sigma(s) + \sigma^p\,.
\label{eq:shiftedp}
\end{equation} Then using H\"older's inequality, \eqref{eq:shiftedp}, \eqref{eq:pq=} and the assumption that $\|D {\bf w}  \|_{L^{p}(Q_r;\R^{Nn})}\le 1$, we obtain
$$\begin{aligned}
r^\alpha\left(\fiint_{Q_r} |D{\bf w}|^{q_*}\, \mathrm{d}z\right)^{\frac{q}{q_*}}
& \le r^\alpha\left(\fiint_{Q_r} |D{\bf w}|^{p\theta}\, \mathrm{d}z\right)^{\frac{1}{\theta}} \left(\fiint_{Q_r} |D{\bf w}|^{p}\, \mathrm{d}z\right)^{\frac{q-p}{p}}\\
&\le r^{\alpha-\frac{(q-p)(n+2)}{p}}\left(\fiint_{Q_r} |D{\bf w}|^{p\theta}\, \mathrm{d}z\right)^{\frac{1}{\theta}}\\
&\le c\left(\fiint_{Q_r} H^-_{\sigma}(|D{\bf w}|)^{\theta}\, \mathrm{d}z\right)^{\frac{1}{\theta}}+ c r^{\alpha-\frac{(q-p)(n+2)}{p}}\sigma^p\,.
\end{aligned}
$$
This completes the proof.
\end{proof}

\section{Caccioppoli inequality via mollification} \label{sec:caccioppoli}

In this section, we derive a Caccioppoli type estimate for $\bfu_{\bm\ell}=\bfu -\bm\ell$ with the shifted function $H_{1+|D\bm\ell|}(z,s)$, where $\bfu$ is the weak solution to \eqref{system1} and $\bm\ell$ is any linear function \eqref{def_l}. We emphasize that in the definition of weak solution, we only assume that $H(\cdot, |D\bfu |)\in L^{1}_{\loc}(\Omega_T)$. Therefore, a more general regularization argument than the Steklov average is needed. To this end, we use mollification in $\R^{n+1}$ in the parabolic setting. We refer to \cite[Chapter 4]{HH} and \cite{HasOk25}.

Let $\tilde\varrho_m\in C^\infty_0(\R^m)$ be a standard mollifier in $\R^m$ defined by $\tilde\varrho_m(y)= \gamma_m \exp(-\tfrac{1}{1-|y|^2})$ if $|y|< 1$ and $\tilde\varrho_m(y)= 0$ if $|y|\ge 1$, where the constant $\gamma_m>0$ is chosen such that $\|\tilde\varrho_m\|_{1}=1$. Then we set $\varrho(z)=\varrho(x,t):=\tilde\varrho_n(x)\tilde\varrho_1(t)$ and, for $h>0$, $\varrho_h(z)=\varrho_h(x,t):= h^{-n-2}\tilde\varrho_n(h^{-1} x)\tilde\varrho_1(h^{-2}t)$.
For $f\in L^1_{\loc}(\Omega_T)$,  $z\in \Omega_T$ and $h>0$ such that $Q_h(z)\Subset \Omega_T$, we define 
\begin{align*}
	f_h(z) := \iint_{Q_h(0)} f(z- \tilde z  ) \varrho_h (\tilde z) \,\d \tilde z\,.
\end{align*}
 
Let $H$ be defined as in \eqref{eq:H}. We observe that the conditions 
\eqref{eq:holdercond} and \eqref{eq:pq=} imply that for both $H$ and $H^*$ holds the parabolic version 
of the (A1) condition in \cite[Definition~4.1.1]{HH}, where the ball $B$ 
is replaced by the parabolic cylinder $Q_r$ in the definition of the 
usual (A1) condition; see \cite[Proposition~7.2.2 and Lemma~4.1.7]{HH}.  
Therefore, the convergence of mollification in the spaces $L^H$ and $L^{H^*}$  in the parabolic setting follows as in \cite[Theorem~4.4.7]{HH}. 
We also refer to \cite[Proposition~4.1]{HasOk25} and the discussion therein for further details. 
 Precisely, if the function $H$ defined in \eqref{eq:H} satisfies 
\eqref{eq:holdercond} and \eqref{eq:pq=}, then for every 
$f\in L^H_{\loc}(\Omega_T)$ and $g\in L^{H^*}_{\loc}(\Omega_T)$ we have
\begin{equation}\label{convergence}
f_h \longrightarrow f 
\quad \text{in }\ L^H_{\loc}(\Omega_T)
\quad \text{and}\quad
g_h \longrightarrow g  
\quad \text{in }\ L^{H^*}_{\loc}(\Omega_T), 
\qquad \text{as }\ h \to 0^+. 
\end{equation}

\begin{lemma}\label{sec3:lem:1}
Let $\bfu$ be a weak solution to \eqref{system1} satisfying \eqref{eq:1.8ok1}, \eqref{eq:1.8ok2} and \eqref{eq:Acontinuity} with \eqref{eq:H}, \eqref{eq:holdercond} and \eqref{eq:pq=}. For every pair of concentric cylinders
$Q_{r_1,\rho_1}(z_0) \subset Q_{r_2,\rho_2}(z_0) \Subset\Omega_T$ with $z_0=(x_0,t_0)$, $0<r_1<r_2\le 1$ and $0<\rho_1<\rho_2 \le1$,  and linear map $\bm\ell$, 
 we have 
	\begin{equation}\label{eq:caccio}\begin{aligned}
&\underset{t\in (t_0-\rho_1, t_0+\rho_1)}{\mathrm{ess\,sup}} \int_{B_{r_1}(x_0)}  {|\bfu_{\bm \ell}(x,t)|^2}\,\mathrm{d}x+\iint_{Q_{r_1,\rho_1}(z_0)} H_{1+ |D{\bm \ell}|}(z,|D \bfu_{\bm \ell}|)\,\mathrm{d}z \\
&\qquad\le \hat{c} \iint_{Q_{r_2,\rho_2}(z_0)} \left[\frac{|\bfu_{\bm \ell}|^2}{\rho_2-\rho_1} + H_{1+ |D {\bm \ell}|}\left(z,\frac{|\bfu_{\bm \ell}|}{r_2-r_1}\right) \right]\,\mathrm{d}z\\
&\qquad\qquad  +\hat{c} |Q_{r_2,\rho_2}| 
\max\Big\{C_1^{\frac{\tilde p}{\tilde p-1}},C_1^{\frac{\tilde q}{\tilde q-1}}\Big\}
H(z^+,1+|D\bm\ell|) \,.
	\end{aligned}\end{equation}
for some $\hat{c}=\hat{c}(n,N,p,q,L,\nu)>0$,
where 
\begin{equation}\label{CrtauP}
C_1 := (r_2+ \sqrt{\rho_2})^{\beta_0} + (r_2+ \sqrt{\rho_2})^{\alpha} (1+|D\bm\ell|)^{q-p} 
\end{equation}
and $1< \tilde p \le \tilde q $ are given in \eqref{phiapq}.
\end{lemma}


\begin{proof}
Let $h_0\in (0,1)$ sufficiently small so that  
$$
h\le  \min\left\{  \frac{1}{4}\mathrm{dist}(Q_{r_2,\rho_2}(z_0),\partial(\Omega_T)), \frac{\rho_2-\rho_1}{4}\right\},
$$
and choose $h\in (0, h_0)$. Let $z^-\in \overline{Q_{r_2,\rho_2}(z_0)}$ such that $a(z^-)=\inf_{z\in Q_{r_2,\rho_2}(z_0) }a(z)$. Since  $\partial_t \bm\ell= \div \bA(z^-, D\bm\ell) =\zero$, we obtain from the weak formulation \eqref{eq:weakformul} that
\begin{equation}\label{sec3:2}
-\iint_{Q_{r_2,\rho_2}(z_0)} [\bfu_{\bm\ell}]_h \cdot \partial_t {\bm\zeta} \,\d z + \iint_{Q_{r_2,\rho_2}(z_0)}   [\bA(\cdot , D\bfu)- \bA(z^- , D\bfu)]_h : D{\bm\zeta} \, \d z =0 
\end{equation}
for all ${\bm\zeta}\in C^\infty_{\mathrm c}(Q_{r_2,\rho_2}(z_0),\R^N)$, where $[\,\cdot\,]_h$ is the mollification stated in above. Note that by approximation the test function ${\bm\zeta}$ can be chosen as a Lipschitz continuous function in $Q_{r_2,\rho_2}(z_0)$ with zero value on $\partial(Q_{r_2,\rho_2}(z_0))$, and $\bA(\cdot , D\bfu), \bA(z^- , D\bfu) \in L^{H^*}_{\loc}(\Omega_T,\R^{N\times n})$.

	Let $\xi\in C_0^\infty(B_{r_2}(x_0))$ be a cut-off function with
	\begin{align}\label{sec3:1}
		0\le \xi\le 1,\quad\xi\equiv1\text{ in } B_{r_1}(x_0)\quad\text{and}\quad\lVert D \xi\rVert_{L^\infty}\le \frac{2}{r_2-r_1},
	\end{align}
and $\eta\in C_0^\infty(I_{\rho_2-h_0}(t_0))$ with
	\begin{align}\label{sec3:3}
			0\le \eta\le 1,\quad \eta\equiv1\text{ in } I_{\rho_1}(t_0)\quad\text{and}\quad\lVert\eta'\rVert_{L^\infty}\le \frac{3}{\rho_2-\rho_1}.
	\end{align}
In addition, for $t_*\in I_{\rho_1}(t_0)$ and $\delta\in(0,h_0)$, define $\eta_\delta: \R\to \R$ as
\begin{align}\label{sec3:4}
	\eta_{\delta}(t)=
	\begin{cases}
		1,&t\in (-\infty,t_*-\delta),\\
		1-\frac{t-t_*+\delta}{\delta},&t\in[t_*-\delta,t_*],\\
		0,&t\in(t_*,\infty).
	\end{cases}
\end{align}
Then we can choose $\bm\zeta=[\bfu_{\bm \ell}]_h\xi^{\tilde{q}}\eta^2\eta_{\delta}$ as a test function in \eqref{sec3:2}, so that
\begin{align}\label{sec3:6}
	\begin{split}
		0 &=
		\iint_{Q_{r_2,\rho_2}(z_0)}\pa_t[\bfu_{\bm \ell}]_h\cdot\bm\zeta\,\mathrm{d}z
		+\iint_{Q_{r_2,\rho_2}(z_0)}[\bA(\cdot, D\bfu)-\bA(\cdot, D \bm\ell)]_h:D\bm\zeta \,\mathrm{d}z\\
		& \qquad + \iint_{Q_{r_2,\rho_2}(z_0)}[\bA(\cdot, D\bm\ell)-\bA(z^-, D \bm\ell)]_h:D\bm\zeta \,\mathrm{d}z\\
		& =: J_1 + J_2 + J_3\,.
	\end{split}
\end{align}

We start by estimating $J_1$. Integration by parts gives
\begin{align*}
		J_1	&=\iint_{Q_{r_2,\rho_2}(z_0)}\tfrac{1}{2}\pa_t(|[\bfu_{\bm\ell}]_h|^2)\xi^{\tilde{q}}\eta^2\eta_{\delta}\,\mathrm{d}z \\
		&=-\iint_{Q_{r_2,\rho_2}(z_0)}|[\bfu_{\bm\ell}]_h|^2\xi^{\tilde{q}}\eta\eta'\eta_{\delta}\,\mathrm{d}z  + \tfrac{1}{2} \fint_{t_*-\delta}^{t_*}\int_{B_{r_2}(x_0)}|[\bfu_{\bm\ell}]_h|^2\xi^{\tilde{q}}\eta^2\,\mathrm{d}x\, \d t \,.
\end{align*}
Then since $\bfu_{\bm\ell} \in L^2_{\loc}(\Omega_T)$, letting $h\to 0^+$ and then $\delta\to 0^+$ and using \eqref{sec3:3} and \eqref{sec3:4}, we obtain 
\begin{equation}\label{J1}\begin{aligned}
	\lim_{\delta\to0^+}\lim_{h\to0^+}J_1 & = - \int_{t_0-\rho_2}^{t^*} \int_{B_{r_2}(x_0)} |\bfu_{\bm\ell}|^2\xi^{\tilde{q}}\eta\eta'\,\mathrm{d}x \, \d t  + \tfrac{1}{2} \int_{B_{r_2}(x_0)}|\bfu_{\bm\ell}(x,t^*)|^2\xi^{\tilde{q}} \,\mathrm{d}x  \\
	&\ge  - c \iint_{Q_{r_2,\rho_2}(z_0)}\frac{|\bfu_{\bm\ell}|^2}{\rho_2-\rho_1} \,\mathrm{d}z   + \tfrac{1}{2} \int_{B_{r_2}(x_0)}|\bfu_{\bm\ell}(x,t^*)|^2\xi^{\tilde{q}}\,\d x\,.
\end{aligned}\end{equation}
Now we estimate $J_2$. Since $\bA(\cdot, D\bfu)-\bA(\cdot, D \bm\ell)\in L^{H^*}_{\loc}(\Omega_T)$ and $D\bfu_{\bm\ell}, \bfu_{\bm\ell} \in L^{H}_{\loc}(\Omega_T)$, by \eqref{convergence},
\begin{align}\label{sec3:11}
	\begin{split}
		\lim_{\delta\to0^+}\lim_{h\to0^+}J_2
		&=\lim_{\delta\to0^+}\lim_{h\to0^+} \bigg[\iint_{Q_{r_2,\rho_2}(z_0)}[\bA(\cdot, D\bfu)-\bA(\cdot, D \bm\ell)]_h:[ D \bfu_{\bm\ell}]_h\xi^{\tilde{q}}\eta^2\eta_\delta\,\d z\\
		&\  \qquad\qquad\qquad+\tilde{q}\iint_{Q_{r_2,\rho_2}(z_0)}[\bA(\cdot, D\bfu)-\bA(\cdot, D \bm\ell)]_h:D\xi \otimes [\bfu_{\bm\ell}]_h\xi^{\tilde{q}-1}\eta^2\eta_\delta\,\d z \bigg]\\
		& 
= \int_{t_0-\rho_2}^{t^*} \int_{B_{r_2}(x_0)}[\bA(z, D\bfu)-\bA(z, D \bm\ell)]: D \bfu_{\bm\ell}\,\xi^{\tilde{q}}\eta^2\, \d x\, \d t \\
&\qquad\qquad +\tilde{q}  \int_{t_0-\rho_2}^{t^*} \int_{B_{r_2}(x_0)} [\bA(z, D\bfu)-\bA(z, D \bm\ell)] :D\xi \otimes \bfu_{\bm\ell}\,\xi^{\tilde q-1}\eta^2\,\d x \, \d t\,.
	\end{split}
\end{align}
Applying \eqref{eq:1.9ok} with $\bm\xi_1=D\bfu$, $\bm\xi_2=D \bm\ell$ and \eqref{(2.6b)} to get
\begin{align*}
	\begin{split}
			&\int_{t_0-\rho_2}^{t^*} \int_{B_{r_2}(x_0)}[\bA(z, D\bfu)-\bA(z, D \bm\ell)]: D \bfu_{\bm\ell}\xi^{\tilde{q}}\eta^2\,\d x \,\d t 	\\
			&\qquad\ge {\tilde{\nu}}\int_{t_0-\rho_2}^{t^*}\int_{B_{r_2}(x_0)}H''(z, 1+|D \bfu|+|D \bm\ell|)|D \bfu_{\bm\ell}|^2\xi^{\tilde{q}}\eta^2\,\d x\,\d t \\
&\qquad\ge {\tilde{c}} \int_{t_0-\rho_2}^{t_*}\int_{B_{r_2}(x_0)}H_{1+|D \bm\ell|}(z, |D \bfu_{\bm\ell}|)\xi^{\tilde{q}}\eta^2\,\d x\,\d t 
	\end{split}
\end{align*}
for some small $\tilde c>0$. To estimate the second term on the right hand side in \eqref{sec3:11}, we use \eqref{eq:(3.4)} and \eqref{sec3:1}, Young's inequality applied to $\varphi(s)=H_{1+|D \bm\ell|}(z,s)$ for any fixed $z$ with \eqref{eq:hok2.4} to conclude that for every $\epsilon>0$,
\begin{align*}
	\begin{split}
		&\int_{t_0-\rho_2}^{t_*}\int_{B_{r_2}(x_0)} [\bA(z, D\bfu)-\bA(z, D \bm\ell)]:D\xi \otimes \bfu_{\bm\ell} \xi^{\tilde{q}-1}\eta^2\,\d x\,\d t\\
		&\quad\ge - c \int_{t_0-\rho_2}^{t_*}\int_{B_{r_2}(x_0)}  H'_{1+|D \bm\ell|}(z,|D \bfu_{\bm \ell}|)\frac{|\bfu_{\bm\ell}|}{r_2-r_1}\xi^{\tilde{q}-1} \eta^2\,\d x\,\d t \\
		&\quad \ge -\epsilon \int_{t_0-\rho_2}^{t_*}\int_{B_{r_2}(x_0)} H_{1+|D \bm\ell|}(z,|D \bfu_{\bm \ell}|)\xi^{\tilde{q}}\eta^2\,\d x\,\d t    -c(\epsilon)\iint_{Q_{r_2,\rho_2}(z_0)}H_{1+|D \bm\ell|}\left(z,\frac{|\bfu_{\bm\ell}|}{r_2-r_1}\right)\eta^2\,\d z\,. 
	\end{split}
\end{align*}
Therefore, it follows that
\begin{align}\label{J2}
	\begin{split}
\lim_{\delta\to0^+}\lim_{h\to0^+}J_2 & \ge \frac{\tilde c}{2} \int_{t_0-\rho_2}^{t_*}\int_{B_{r_2}(x_0)}  H_{1+|D \bm\ell|}(z,|D \bfu_{\bm \ell}|)\xi^{\tilde{q}}\eta^2\,\d x\,\d t \\
&\qquad -c\iint_{Q_{r_2,\rho_2}(z_0)}H_{1+|D \bm\ell|}\left(z,\frac{|\bfu_{\bm\ell}|}{r_2-r_1}\right)\,\d z.
	\end{split}
\end{align}
We finally estimate $J_3$. Similarly to the estimation of $J_2$, we have 
\begin{align*}
		\lim_{\delta\to0^+}\lim_{h\to0^+} J_3 & =\int_{t_0-\rho_2}^{t^*} \int_{B_{r_2}(x_0)}[\bA(z, D\bm\ell)-\bA(z^-, D \bm\ell)]: [D \bfu_{\bm\ell}\xi^{\tilde{q}}+\tilde{q} D\xi \otimes \bfu_{\bm\ell}\xi^{\tilde{q}-1}]\eta^2\, \mathrm{d} x\, \mathrm{d} t \,.
\end{align*}
By \eqref{eq:Acontinuity} and \eqref{eq:holdercond}, for $z\in Q_{r_2,\rho_2}(z_0)$, 
\begin{align*}
|\bA(z, D\bm\ell)-\bA(z^-, D \bm\ell)| 
& \le c (r_2 + \sqrt{\rho_2})^{\beta_0} H'(z,1+|D\bm\ell|) +  c (r_2 + \sqrt{\rho_2})^{\alpha} (1+|D\bm\ell|)^{q-1}  \\
& \le c C_1  H_{1+|D\bm\ell|}'(z,1+|D\bm\ell|),
\end{align*}
where $C_1$ is given in \eqref{CrtauP}.
Plugging this in the previous identity and using Young's inequality, we obtain
\begin{equation}\label{sec3:9}
\begin{aligned}
\lim_{\delta\to0^+}\lim_{h\to0^+} J_3  & \ge - \frac{\tilde c}{4} \int^{t^*}_{t_0-\rho_2}  \int_{B_{r_2}(x_0)} H_{1+|D\bm\ell|} (z, |D\bfu_{\bm\ell}|)\xi^{\tilde{q}}\eta^2\, \mathrm{d} x\, \mathrm{d} t \\
& \qquad -c\iint_{Q_{r_2,\rho_2}(z_0)}H_{1+|D \bm\ell|}\left(z,\frac{|\bfu_{\bm\ell}|}{r_2-r_1}\right)\,\d z\\
& \qquad - c |Q_{r_2,\rho_2}|\max\Big\{C_1^{\frac{\tilde p}{\tilde p-1}},C_1^{\frac{\tilde q}{\tilde q-1}}\Big\}  H(z^+,1+|D\bm\ell|) \,.
\end{aligned}\end{equation}

Therefore, combining \eqref{sec3:6}, \eqref{J1}, \eqref{J2} and  \eqref{sec3:9} yields
$$\begin{aligned}
\int_{B_{r_2}(x_0)}  |\bfu_{\bm\ell}(x,t^*)|^2\xi^{\tilde{q}}\,\d x  + & \int^{t^*}_{t_0-\rho_2}  \int_{B_{r_2}(x_0)}   H_{1+|D\bm\ell|} (z, |D\bfu_{\bm\ell}|)  \xi^{\tilde{q}}\eta^2\, \mathrm{d} x\, \mathrm{d} t       \\
& \le  c\iint_{Q_{r_2,\rho_2}(z_0)}\left[\frac{|\bfu_{\bm\ell}|^2}{\rho_2-\rho_1} + H_{1+|D \bm\ell|}\left(z,\frac{|\bfu_{\bm\ell}|}{r_2-r_1}\right)\right]\,\d z\\
& \qquad + c |Q_{r_2,\rho_2}|\max\Big\{C_1^{\frac{\tilde p}{\tilde p-1}},C_1^{\frac{\tilde q}{\tilde q-1}}\Big\}   H(z^+,1+|D\bm\ell|) \,.
\end{aligned}$$
Since $t_*\in (t_0-\rho_1,t_0+\rho_1)$ is arbitrary, considering  \eqref{sec3:1} and \eqref{sec3:3}, we obtain \eqref{eq:caccio}. 
\end{proof}



\section{$\mathcal{A}$-caloric approximation and linearization}\label{sec:approxandlin}

The $\mathcal A$-caloric approximation, as stated in \cite[Theorem 3.8]{OkSciStroffo2024}, is a crucial tool in the nondegenerate setting. It will be used in Section~\ref{sec:linearization} to establish a comparison estimate between the gradient of a weak solution to the system \eqref{system1} and the gradient of a smooth 
$\mathcal A$-caloric map (see Section~\ref{subsec:CZ}). This comparison is possible provided that 
$\bfu$ can be shown to be an almost solution of a suitable 
$\mathcal A$-caloric system.

\subsection{Regularity estimates for linear systems with constant coefficients}\label{subsec:CZ}
We first introduce an $\A$-caloric map and it regularity estimates.
Let $ \A  = (\A^{\alpha\beta}_{ij}) \in \R^{N^2n^2}$ satisfy the \emph{Legendre-Hadamard condition}: for every $\mathbf{a}= (a^1,\dots,a^N)\in \R^N$ and  $\mathbf{b}=(b_1,\dots,b_n)\in \R^n$,
\begin{equation*}\label{constantA}
\A ({\bf a} \otimes {\bf b})  \, : \, ({\bf a}\otimes {\bf b})   =\A^{\alpha\beta}_{ij} a^\alpha a^\beta b_i  b_j \ge \nu |\mathbf{a}|^2 |\mathbf{b}|^2 
\end{equation*}
for some $\nu>0$.  Then a weak solution $\bfh : Q_r \to \R^N$ to the linear system with coefficient $\A$, given by
$$
\bfh_t -\div (\A D\bfh) = {\bf 0} \quad \text{in }\ Q_r\,,
$$
is called an \textit{$\A$-caloric map}. 
By standard regularity theory (see \cite{Cam66,DiB_book}), $\bfh \in C^\infty(Q_r,\R^N)$. Moreover, one can derive the following estimates (see, e.g., \cite[(5.9)--(5.12)]{Cam66,DiB_book}), which will be used later.
\begin{lemma}\label{Lem:Acaloricmap}
Suppose $\bfh\in C^\infty(Q_r,\R^N)$ is an $\A$-caloric map in $Q_r$. Then we have that 
\begin{equation}\label{supestimate}
\sup_{Q_{r/2}} \left(|\bfh|+r |D\bfh|+r^2|D^2\bfh|\right)  \le c  \fiint_{Q_{r}} |\bfh|  \, \d z\,,
\end{equation}
and
\begin{equation}\label{DhLip}
\sup_{Q_{r/2}} |D\bfh|  \le c  \fiint_{Q_{r}} |D\bfh|  \, \d z\,.
\end{equation}
Here, the constants $c>0$ depend on $n, N, \nu$ and $|\mathcal A|$.
\end{lemma}

We recall the following $\A$-caloric approximation result, proved in \cite[Theorem 3.8]{OkSciStroffo2024}. 

\begin{theorem}\label{thm:Acaloric}  ($\mathcal A$-caloric approximation)
Let $\mu,\sigma_0,C_0>0$ and $\psi$ be an $N$-function with $\Delta_2(\psi,\psi^*)<\infty$. For every $\epsilon\in(0,1)$, there exists $\delta>0$ depending on $\sigma_0$, $C_0$, $\Delta_2(\psi,\psi^*)$ and $\epsilon$   such that
if $\bfu\in L^1(-r^2,r^2,W^{1,1}(B_r;\R^N))$ with $\psi^{1+\sigma_0}(|D\bfu|)\in L^{1}(Q_r)$ and $\psi^{1+\sigma_0} ({\bf H}) \in L^{1}(Q_r)$ satisfy
$$
\partial_t \bfu = \div {\bf H}  \quad \text{in }\ Q_r,
$$
in the distributional sense, and the inequality
$$
\left( \fiint_{Q_r} \left[ \psi(|D\bfu|)^{1+\sigma_0}+ \psi(|{\bf H}|)^{1+\sigma_0}\right]\, \d z \right)^{\frac{1}{1+\sigma_0}}  \le C_0 \psi(\mu),
$$
and for every $\bm\zeta\in C^\infty(\overline{Q_r};\R^{N})$ with $\bm\zeta={\bf 0}$ on $\partial B_r \times (-r^2,r^2)$,
$$
\frac{1}{|Q_r|}\left|  \iint_{Q_r} \bfu \cdot \bm\zeta_t - [\mathcal A D\bfu : D\bm\zeta] \,\d z - \left[\int_{B_r} \bfu \cdot \bm\zeta_t \, \d x\right]^{t=r^2}_{t=-r^2}\right| \le \delta \mu  \|D\bm\zeta\|_{L^\infty(Q_r,\R^{N\times n})},
$$
then 
$$
\fint_{Q_r} \psi(|D\bfu-D\bfh|)\,\d z \le \epsilon \psi(\mu),
$$
where $\bf h$ is the weak solution to
$$
\begin{cases}
\partial_t \bfh- \div (\A D\bfh) = 0  \quad \text{in}\ \ Q_r, \\
 \bfh =\bfu \quad \text{on}\ \ \partial_{\mathrm{p}} Q_r.
\end{cases}$$
\end{theorem}

\subsection{Linearization} \label{sec:linearization}

To take advantage of the comparison estimates in Theorem \ref{thm:Acaloric}, we demonstrate that the weak solution 
$\bfu$ to \eqref{system1} is an 'almost' weak solution of a linear system with constant coefficients.
We first note that, since $\|H(\cdot, |D\bfu|)\|_{L^1(\Omega_T)} < \infty$ by the definition of the weak solution to \eqref{system1}, there exists small $R_0\in (0,1)$ such that 
\begin{equation}
\iint_{Q_r(z_0)}H(z,1+|D\bfu|)\, \d z \le 1
\quad \text{for every  } \ Q_r(z_0)\Subset \Omega_T \ \text{ with } \ r \in (0,R_0]\,.
\label{eq:smallnesshyp}
\end{equation}

\begin{lemma}\label{lem:linearization}
Let $\bfu$ be a weak solution to \eqref{system1} satisfying \eqref{eq:1.8ok1}, \eqref{eq:1.8ok2}, \eqref{eq:Acontinuity} and \eqref{offdiagonal} with \eqref{eq:H}, \eqref{eq:holdercond} and \eqref{eq:pq}, $\bm\ell$ be any linear map defined in  \eqref{def_l}, and $Q_r=Q_r(z_0) \Subset \Omega_T$ with $z_0=(x_0,t_0)$ and $r\in (0,R_0]$.
Then for every $\bm\zeta \in C^\infty(\overline{Q_r};\R^N)$ with $\bm\zeta= {\bf 0} $ on $\partial B_r \times I_r$ and $\|D\bm\zeta\|_{L^\infty(Q_r;\R^{N\times n})} \le 1$, we have 
\begin{equation}\label{Acalori_pf1}
\begin{aligned}
&\frac{1}{|Q_r|}\left|\iint_{Q_r} \bfu_{\bm\ell}  \cdot \bm\zeta_t - \langle D_{\bm\xi}\bA(z^-,D\bm\ell)  D\bfu_{\bm\ell} \,|\, D\bm\zeta\rangle \, \d z  - \left[\int_{B_r} \bfu_{\bm\ell}  \cdot \bm\zeta  \, \d x \right]^{t=t_0+r^2}_{t=t_0-r^2}  \right| \\
& \le  c r^{\beta_1}  \left\{ \big((H^-)'\circ (H^{-})^{-1}\big)\left(\fiint_{Q_r}H_{1+|D\bm\ell|}^-(|D\bfu_{\bm\ell}|)\,\mathrm{d}z\right) + (H^-)'(1+|D\bm\ell|)\right\}\\
&\qquad + c (H^-)'(1+|D\bm\ell|)\left\{ \fiint_{Q_r} \frac{H^-_{1+|D{\bm\ell}|}(|D\bfu_{\bm\ell}|)}{H^-(1+|D\bm\ell|)}\, \d z+  \left(  \fiint_{Q_r}  \frac{H^-_{1+|D\bm\ell|}(|D\bfu_{\bm\ell}|)}{H^-(1+|D{\bm\ell}|)} \, \d z \right)^{\frac{1+\gamma}{2}} \right\}
\end{aligned}\end{equation}
for some $c=c(n,N,p,q,L,\nu)$, where $\beta_1=\beta_1(n,p,q,\alpha,\beta_0)>0$ is defined as
\begin{equation}
\beta_1:=\min \left\{\beta_0,\alpha-(n+2)\frac{q-p}{p}\right\}\,,
\label{def:beta1}
\end{equation}
and $\beta_0$, $\gamma$,  $\bfu_{\bm\ell}$ and $H^\pm :=H_{Q_r}^\pm$ are from \eqref{eq:Acontinuity}, \eqref{offdiagonal}, \eqref{def_ul} and \eqref{eq:Hpmdef}, respectively. 
\end{lemma}

\begin{proof} Throughout the proof, we will denote $z^-_{z_0,r}$ by $z^-$. 
Note that by a similar approximation used in the proof of Lemma~\ref{sec3:lem:1}, we have from \eqref{eq:weakformul} that
\[
-\iint_{Q_r} \bfu \cdot {\bm\zeta}_t \,\d z +  \left[\int_{B_r} \bfu  \cdot \bm\zeta  \, \d x \right]^{t=t_0+r^2}_{t=t_0-r^2} + \iint_{Q_r}   \bA(z , D\bfu) : D{\bm\zeta} \, \d z =0\,.
\]
From this and the fact that ${\bm\ell}_t=\div (\bA(z^-,D{\bm\ell}))=\zero$, we obtain
\begin{equation}
\begin{aligned}
&\frac{1}{|Q_r|}\left|\iint_{Q_r} \bfu_{\bm\ell}  \cdot \bm\zeta_t - \langle D_{\bm\xi}\bA (z^-,D{\bm\ell})   D\bfu_{\bm\ell} \, |\,   D\bm\zeta\rangle \, \d z  - \left[\int_{B_r} \bfu_{\bm\ell}  \cdot \bm\zeta   \, \d x \right]^{t=r^2}_{t=-r^2}\right|\\
&=\frac{1}{|Q_r|}\left| \iint_{Q_r} \bfu  \cdot \bm\zeta_t - \langle D_{\bm\xi}\bA (z^-,D{\bm\ell})  D\bfu_{\bm\ell} \,|\, D\bm\zeta\rangle \, \d z  - \left[\int_{B_r} \bfu \cdot \bm\zeta  \, \d x \right]^{t=r^2}_{t=-r^2}\right| \\
&= \left|\fiint_{Q_r} \langle \bA(z,D\bfu) -\bA(z^-,D{\bm\ell}) \,|\, D\bm\zeta \rangle -  \langle D_{\bm\xi}\bA (z^-,D{\bm\ell})  D\bfu_{\bm\ell} \,|\, D\bm\zeta\rangle \, \d z \right|\\
&=\left|\fiint_{Q_r} \langle \bA(z,D\bfu) -\bA(z^-,D
\bfu) \,|\, D\bm\zeta \rangle  \, \d z \right.\\
&\qquad \left. +\fiint_{Q_r} \int_0^1 \langle [D_{\bm\xi}\bA( z^-,sD\bfu_{\bm\ell}+ D{\bm\ell}) - D_{\bm\xi}\bA (z^-,D{\bm\ell})] D\bfu_{\bm\ell} , D\bm\zeta\rangle \, \d s \, \d z \right|\\
&\leq \fiint_{Q_r} | \bA(z,D\bfu) -\bA(z^-,D
\bfu) |  \, \d z\\
&\qquad + \fiint_{Q_r} \left[\int_0^1 | D_{\bm\xi}\bA( z^-,sD\bfu_{\bm\ell}+ D{\bm\ell} ) - D_{\bm\xi}\bA (z^-,D{\bm\ell})|  \, \d s \right]|D\bfu_{\bm\ell} | \, \d z\\
&=: I_1+I_2 \,.
\end{aligned}
\label{eq:estimalmostcal}
\end{equation}

We first estimate $I_1$. By \eqref{eq:Acontinuity} and the definition of $z^-$,
\begin{equation}\label{eq:I1}\begin{aligned}
I_1 & \le c r^{\beta_0} \fiint_{Q_r} H'(z,1+|D\bfu|)\, \mathrm{d}z + c\fiint_{Q_r} (a(z)-a(z^-))(1+|D\bfu|)^{q-1}\,\mathrm{d}z\\
& \le c r^{\beta_0} \fiint_{Q_r} (H^-)'(1+|D\bfu|)\, \mathrm{d}z + cr^\alpha \fiint_{Q_r} (1+|D\bfu|)^{q-1}\,\mathrm{d}z. 
\end{aligned}\end{equation}
By Jensen's inequality and  the fact that $((H^-)^*)^{-1} \approx  (H^-)'\circ (H^-)^{-1}$, the first integral is estimated as 
$$\begin{aligned}
\fiint_{Q_r} (H^-)'(1+|D\bfu|)\, \mathrm{d}z 
&\le  ((H^-)^*)^{-1}\left( \fiint_{Q_r} \big((H^-)^*\circ (H^-)'\big)(1+|D\bfu|)\, \mathrm{d}z \right)\\
& \sim \big((H^-)'\circ (H^-)^{-1}\big)\left( \fiint_{Q_r} H^-(1+|D\bfu|)\, \mathrm{d}z \right)\,.
\end{aligned}$$
For  the second integral, by \eqref{eq:holdercond}, \eqref{eq:pq}, \eqref{eq:smallnesshyp} and Jensen's inequality for the function $t\mapsto H^-(t^{1/p})$,
$$\begin{aligned}
r^{\alpha}\fiint_{Q_r}(1+|D\bfu|)^{q-1}\,\mathrm{d}z &\le r^{\alpha}\left(\fiint_{Q_r}(1+|D\bfu|)^{p}\,\mathrm{d}z\right)^{\frac{q-1}{p}}\\
&\le r^{\alpha-\frac{q-p}{p}(n+2)} \left(\fiint_{Q_r}(1+|D\bfu|)^{p}\,\mathrm{d}z\right)^{\frac{p-1}{p}}\\
&\le r^{\alpha-\frac{q-p}{p}(n+2)} (H^{-})^{-1}\left(\fiint_{Q_r}H^-(1+|D\bfu|)\,\mathrm{d}z\right)^{p-1}\\
&\le r^{\alpha-\frac{q-p}{p}(n+2)} \big((H^-)'\circ (H^{-})^{-1}\big)\left(\fiint_{Q_r}H^-(1+|D\bfu|)\,\mathrm{d}z\right)\,.
\end{aligned}$$
Therefore, inserting the above estimates into \eqref{eq:I1} and  using \eqref{eq:approx} for $\varphi(t):=H^-(t)$, $\sigma:= 1+ |D \bm\ell|$ and $s=|D\bfu_{\bm\ell}|$, we have 
\begin{equation}\label{eq:estimI1}\begin{aligned}
I_1& \le c r^{\beta_1}   \big((H^-)'\circ (H^{-})^{-1}\big)\left(\fiint_{Q_r}H^-(1+|D\bfu|)\,\mathrm{d}z\right) \\
& \le c r^{\beta_1}   \big((H^-)'\circ (H^{-})^{-1}\big)\left(\fiint_{Q_r}H^-_{1+|D\bm\ell|}(|D\bfu_{\bm\ell}|)\,\mathrm{d}z + H^-(1+|D\bm\ell|)\right)\,,
\end{aligned}\end{equation}
where $\beta_1>0$ is defined in \eqref{def:beta1}. 

We next estimate $I_2$ by decomposing $Q_r$ into
 $F_1= \{z\in Q_r : 2|D\bfu_{\bm\ell}(z)| >1+|D{\bm\ell}|\}$ and  $F_2= \{z\in Q_r : 2|D\bfu_{\bm\ell}(z)|  \le 1+ |D{\bm\ell}|\}$, 
 together with their characteristic functions $\chi_{F_1}$ and $\chi_{F_2}$.
For $z\in F_1$, the inequality
\begin{equation}
|D\bfu(z)|+1+|D{\bm\ell}|\le |D\bfu_{\bm\ell}(z)| + 2(1+|D{\bm\ell}|) \le 5 |D\bfu_{\bm\ell}(z)|
\label{eq:inequbounds}
\end{equation}
holds. Using this along with \eqref{eq:1.8ok1}, \eqref{DE08lem20} with $\phi(t)=H_1^-(t)$,  \eqref{(2.6a)}, \eqref{(2.6b)} and \eqref{eq:inequbounds},
we can estimate:
$$\begin{aligned}
&\int_0^1 | D_{\bm\xi}\bA( z^-,sD\bfu_{\bm\ell}(z) + D{\bm\ell}) - D_{\bm\xi}\bA (z^-,D{\bm\ell})|  \, \d s \\
& \le  c \int_0^1 (H^-)''(1+ |sD\bfu(z)+(1-s)D{\bm\ell}|) \, \d s + c (H^-)'' (1+|D{\bm\ell}|) \\
& \le c \int_0^1 \frac{(H_1^-)'(|sD\bfu(z)+(1-s)D{\bm\ell}|)}{|sD\bfu(z)+(1-s)D{\bm\ell}|} \, \d s + c (H^-)'' (1+|D{\bm\ell}|)  \\
&\le c  \frac{(H^-)'(  1+|D\bfu(z)|+|D{\bm\ell}|)}{1+ |D\bfu(z)| +|D{\bm\ell}|} +c  \frac{(H^-)'( 1+|D\bfu(z)| +|D{\bm\ell}|)}{ 1+|D{\bm\ell}|}  \\
&\le \frac{c}{1+|D{\bm\ell}|}  (H^-)'( |D\bfu_{\bm\ell}(z)| )\\
&\le \frac{c}{1+|D{\bm\ell}|}     (H^-)'( |D\bfu_{\bm\ell}(z)| +1+|D{\bm\ell} | ) \frac{5|D\bfu_{\bm\ell}(z)|}{|D\bfu_{\bm\ell}(z)| + 1+|D{\bm\ell} |} \le  \frac{c}{1+|D{\bm\ell}|}  (H^-_{1+|D{\bm \ell}|})'( |D\bfu_{\bm\ell}(z)|)\,.
\end{aligned}$$
Hence, by \eqref{(2.6a)}, 
\begin{equation}
\begin{aligned}
&\fiint_{Q_r} \chi_{F_1}(z) \left[\int_0^1 | D_{\bm\xi}\bA(z^-, sD\bfu_{\bm\ell}+ D{\bm\ell} ) - D_{\bm\xi}\bA (z^-, D{\bm\ell})|  \, \d s \right]|D\bfu_{\bm\ell}|  \, \d z \\
&\hspace{6cm} \le   \frac{c}{1+|D\bm\ell|} \fiint_{Q_r} H^-_{1+|D{\bm\ell}|}(|D\bfu_{\bm\ell}|)\, \d z \,.
\end{aligned}
\label{eq:estimlargegrad}
\end{equation}
On the other hand, if $z\in F_2$,
applying \eqref{offdiagonal} with $z=z^-$, $\bm\xi_1=D{\bm\ell}$ and $\bm\xi_2=sD\bfu_{\bm\ell}(z) + D{\bm\ell}$ we have
$$
\int_0^1 | D_{\bm\xi}\bA( z^-,sD\bfu_{\bm\ell}(z) + D{\bm\ell}) - D_{\bm\xi}\bA (z^-,D{\bm\ell})|  \, \d s \le c \, \left(\frac{|D\bfu_{\bm\ell}(z)|}{1+|D{\bm\ell}|}\right)^{\gamma} (H^-)''(1+|D{\bm\ell}|) \,
$$
and, since $(H^-)'$ is increasing, using also \eqref{(2.6a)} and \eqref{(2.6b)}, we get 
$$\begin{aligned}
\frac{|D\bfu_{\bm\ell}(z)|^2}{(1+|D{\bm\ell}|)^2} &\le \frac{(H^-)'(|D\bfu_{\bm\ell}(z)|+1+|D{\bm\ell}|)}{(H^-)'(1+|D{\bm\ell}|)(1+|D{\bm\ell}|)} \cdot \frac{|D\bfu_{\bm\ell}(z)|^2}{1+|D{\bm\ell}|} \\
&\le \frac{(H^-)'(|D\bfu_{\bm\ell}(z)|+1+|D{\bm\ell}|)}{H^-(1+|D{\bm\ell}|)} \cdot \frac{|D\bfu_{\bm\ell}(z)|^2}{|D\bfu_{\bm\ell}|+ \frac{1}{2}(1+|D{\bm\ell}|)}  \le c \frac{H^-_{1+|D\bm\ell|}(|D\bfu_{\bm\ell}|)}{H^-(1+|D{\bm\ell}|)}\,.
\end{aligned}$$
Then combining these estimates, using the definition of $F_2$, and  applying H\"older's inequality, we have
\begin{equation}
\begin{aligned}
&\fiint_{Q_r}  \chi_{F_2}(z) \left[\int_0^1 | D_{\bm\xi}\bA( z^-,sD\bfu_{\bm\ell}+ D{\bm\ell} ) - D_{\bm\xi}\bA (z^-,D{\bm\ell})|  \, \d s \right] |D\bfu_{\bm\ell} | \, \d z \\
&\le c (H^-)'(1+|D{\bm\ell}|) \fiint_{Q_r}    \chi_{F_2}(z)  \left(\frac{|D\bfu_{\bm\ell}|}{1+|D{\bm\ell}|}\right)^{1+\gamma}\, \d z\\
&\le c (H^-)'(1+|D{\bm\ell}|) \left(  \fiint_{Q_r}  \frac{H^-_{1+|D\bm\ell|}(|D\bfu_{\bm\ell}|)}{H^-(1+|D{\bm\ell}|)} \, \d z \right)^{\frac{1+\gamma}{2}} \,.
\end{aligned}
\label{eq:estimsmallgrad}
\end{equation}
Therefore, collecting \eqref{eq:estimlargegrad} and \eqref{eq:estimsmallgrad}, we obtain
\begin{equation}
I_2\le c (H^-)'(1+|D\bm\ell|)\left\{ \fiint_{Q_r} \frac{H^-_{1+|D{\bm\ell}|}(|D\bfu_{\bm\ell}|)}{H^-(1+|D\bm\ell|)}\, \d z+  \left(  \fiint_{Q_r}  \frac{H^-_{1+|D\bm\ell|}(|D\bfu_{\bm\ell}|)}{H^-(1+|D{\bm\ell}|)} \, \d z \right)^{\frac{1+\gamma}{2}} \right\}\,.
\label{eq:estimI2}
\end{equation}
Finally, plugging \eqref{eq:estimI1} and \eqref{eq:estimI2} into \eqref{eq:estimalmostcal}, we get \eqref{Acalori_pf1}. This concludes the proof.  
\end{proof}

\section{Decay estimate}\label{sec:decayestimate}

For $z_0=(x_0,t_0)\in \Omega_T$ and $\rho>0$ such that $Q_{\rho}(z_0)\Subset\Omega_T$, we define the excess functional 
\begin{equation}
\Phi(z_0,\rho,\bfu) = \fiint_{Q_{\rho}(z_0)} H^-_{1+|(D\bfu)_{z_0,\rho}|}\left( \left|D\bfu -(D\bfu)_{z_0,\rho}\right|\right) \, \d z + \rho^{\frac{\beta_1}{2}}\,,
\label{eq:excess}
\end{equation}
where $H^-_{1+|(D\bfu)_{z_0,\rho}|}:=(H^-_{Q_\rho(z_0)})_{1+|(D\bfu)_{z_0,\rho}|}$ and  $\beta_1$ is given in \eqref{def:beta1},
and linear map
\begin{equation}
  {\bm\ell}_{z_0,\rho,\bfu} (x) := \left(\frac{n+2}{\rho^2}\fiint_{Q_\rho(z_0)} \bfu(x,t) \otimes (x-x_0)\, \d x\, \d t   \right) (x-x_0) + (\bfu)_{Q_{\rho}(z_0)}\,. 
\label{eq:linearmap}
\end{equation}
Note that $\bm\ell_{z_0,\rho,\bfu}$ is the minimizer of 
the functional
$$
F(\bm\ell) = \fiint_{Q_\rho(z_0)} |\bfu- \bm\ell|^2\,\d z\,,
$$
where $\bm\ell=\bm\ell(x)$ is any affine function from $\R^n$ to $\R^N$, and 
\begin{equation}\label{DlDl}
|D\bm\ell_{z_0,\rho,\bfu}-D\bm\ell| = \frac{n+2}{\rho^2}\left|\fiint_{Q_\rho(z_0)}(\bfu - \bm\ell)\otimes (x-x_0) \, \d x\, \d t \right| \le (n+2) \fiint_{Q_\rho(z_0)} \frac{|\bfu - \bm\ell|}{\rho}\,\d z\,,
\end{equation}
see \cite[Section 2.5]{BoDuMin13}.
The following Lemma collects some results of (almost) minimality proven for function $\bm\ell_{z_0,\rho,\bfu}$ in \cite[Lemma 2.6, Lemma 2.7 and Remark 2.6]{FILV23}.
\begin{lemma}
Let $\varphi$ be an $N$-function satisfying $\Delta_2(\varphi,\varphi^*)<\infty$, and let $\bfu\in L^\varphi(Q_\rho(z_0),\mathbb{R}^N)$. There exists a constant $c=c(n,\Delta_2(\varphi,\varphi^*))>0$ such that for every affine function $\bm\ell: \mathbb{R}^n\to \mathbb{R}^N$,
\begin{equation}
\fiint_{Q_\rho(z_0)} \varphi \left(\frac{|\bfu - \bm\ell_{z_0,\rho}|}{\rho}\right)\,\mathrm{d}z \leq c \fiint_{Q_\rho(z_0)} \varphi \left(\frac{|\bfu - \bm\ell|}{\rho}\right)\,\mathrm{d}z\,.
\label{eq:filv1}
\end{equation}
Moreover, if $D\bfu \in L^{\phi}(Q_\rho(z_0),\R^{N\times n})$, then for every ${\bf Q} \in \mathbb{R}^{N\times n}$,
\begin{equation}
\fiint_{Q_\rho(z_0)} \varphi_{1+|(D \bfu)_{z_0,\rho}|}\left({|D \bfu - (D \bfu)_{z_0,\rho}|}\right)\,\mathrm{d}z \leq c \fiint_{Q_\rho(z_0)} \varphi_{1+|{\bf Q}|}\left(|D\bfu - {\bf Q}|\right)\,\mathrm{d}z \,.
\label{eq:filv3}
\end{equation}
\end{lemma}

Now, we prove a decay estimate for the excess $\Phi(z_0,\rho,\bfu)$.

\begin{lemma} \label{Lem:nondegenerate1}
Let $\bfu$ be a weak solution to \eqref{system1}, $K\ge 1$, and $\beta\in(0,\tfrac{\beta_1}{4})$ where $\beta_1$ is given in \eqref{def:beta1}. There exist small $r_0\in (0,R_0]$ and $\delta_0,\tau \in(0,1)$ depending on $n,N,p,q,\alpha,L,\nu,\gamma,K$ and $\beta$ such that if 
$Q_{2r}(z_0)\Subset \Omega_T$ with $r\in(0,r_0]$, 
\begin{equation}
\label{Lem:nondegenerate1_Ass1}
|(D\bfu)_{z_0,r}| \le K
\end{equation}
and
\begin{equation}\label{Lem:nondegenerate1_Ass2}
\Phi(z_0,r,\bfu)  \le \delta_0 \,,
\end{equation}
then
\begin{equation}\label{Dldecay}
|(D\bfu)_{z_0,\tau r}| \le  |(D\bfu)_{z_0,r}| + c_1 \tau^{-n-2} \Phi(z_0, r, \bfu)^{\frac12}
\end{equation}
for some $c_1=c_1(n)>0$,
and
\begin{equation}\label{nondegenerate_decay}
\Phi(z_0,\tau r,\bfu) \le  \tau^{2\beta} \Phi(z_0,r,\bfu)\,.
\end{equation}
\end{lemma}

\begin{proof}
\textit{Step 1.\,(Setting)} All constants $c$, including implicit constants, may depend also on $K$. For simplicity, we write $Q_\rho:=Q_\rho(z_0)$ for any $\rho>0$, $z^-:= z^-_{z_0,r}$,
$H^\pm_\sigma := (H^\pm_{Q_r})_\sigma $, $H^{\pm,*}_\sigma:=(H^{\pm}_\sigma)^*$, $\bm\ell_{\tau r}:=\bm\ell_{z_0,\tau r,\bfu}$ (see \eqref{eq:linearmap}), and  
$$
\tilde{\bm\ell}(x):=(D\bfu)_{z_0,r} (x-x_0) +  (\bfu)_{z_0,r}\,.
$$ 
Note that, by \eqref{Lem:nondegenerate1_Ass1},
\begin{equation}\label{Dlequiv}
1+|D\tilde{\bm\ell}| \sim 1\,.
\end{equation}
Set
\begin{equation}\label{def:mu}
\mu := \left(\fiint_{Q_{r}}H^-_{1+|D\tilde{\bm\ell}|} (|D\bfu_{\tilde{\bm\ell}}|)\,\d z+ r^{\frac{\beta_1}{2}}\right)^{1/2}  =   \Phi(z_0,r,\bfu)^{\frac{1}{2}}\le  \delta_0^{\frac12} \le 1 \,.
\end{equation}
\textit{Step 2.\,($\A$-caloric approximation)} Observe that 
$$
\partial_t \bfu_{\tilde{\bm\ell}} = \partial_t \bfu  = \div {\bf G} \quad \text{in } \ Q_{r}\,, \quad \text{where }\ {\bf G}:=\bA(z,D \bfu) - \bA (z^-,D \tilde{\bm\ell})\,,
$$
in the distributional sense,
 and that
\begin{equation}
\begin{aligned}
| {\bf G}| &\le | \bA(z, D \bfu) - \bA(z^-, D\bfu)| +| \bA(z^-, D \bfu) - \bA(z^-, D \tilde{\bm\ell})|\\
&\le c(a(z)-a(z^-))(1+|D\bfu|)^{q-1}+c r^{\beta_0} H'(z,1+|D\bfu|)
+ (H^-_{1+|D\tilde{\bm\ell} |})'( |D \bfu_{\tilde{\bm\ell}}|)\\
&\le c(a(z)-a(z^-))(1+|D\bfu|)^{q-1}+c r^{\beta_0} (H^-)'(1+|D\bfu|)
+ (H^-_{1+|D\tilde{\bm\ell} |})'( |D \bfu_{\tilde{\bm\ell}}|)\,.
 \end{aligned}
\label{eq:estimG}
\end{equation}
Indeed, the first term above $| \bA(z, D \bfu) - \bA(z^-, D\bfu)|$ can be estimated by \eqref{eq:Acontinuity}, while by using \eqref{eq:(3.4)} with ${\bf P}= D\bfu$ and ${\bf Q}= D\tilde{\bm\ell}$, we obtain
\begin{equation*}
| \bA(z^-, D \bfu) - \bA(z^-, D \tilde{\bm\ell})|  \lesssim (H_{1+|D\tilde{\bm\ell}|}^-)'(|D\bfu_{\tilde{\bm\ell}}|)\,.
\end{equation*}

Let
\begin{equation}\label{def:p0p1}
 p_1 := \min\left\{\frac{p}{q-1},p,\frac{q}{q-1},2\right\}>1
 \quad\text{and}\quad
p_0:=\frac{1+p_1}{2}.
\end{equation}
Then one can see that  both functions  $\psi(s)=H^{-,*}_{\sigma}(s^{1/p_1}) $ and $\psi(s)=H^{-}_{\sigma}(s^{1/p_1})$ with $\sigma\ge 0$ satisfy the assumption in Lemma~\ref{lem:Jensen}.
We first estimate, using \eqref{eq:estimG},  
\begin{equation*}
\begin{split}
\fiint_{Q_r} |{\bf G}|^{p_1}\, \d z  
& \le   c \fiint_{Q_r} \left[(a(z)-a(z^-))(1+|D\bfu|)^{q-1}\right]^{p_1} \, \d z + c \fiint_{Q_r} \left[r^{\beta_0} (H^-)'(1+|D\bfu|)\right]^{p_1}\, \d z\\
&\qquad + c \fiint_{Q_r} (H^-_{1+|D\tilde{\bm\ell} |})'(|D \bfu_{\tilde{\bm\ell}}|)^{p_1}\, \d z \\
&=: II_1+II_2+II_3\,.
\end{split}
\end{equation*} 
By H\"older's inequality, \eqref{eq:holdercond}, \eqref{eq:smallnesshyp}, \eqref{eq:approx}, \eqref{Dlequiv},  \eqref{def:beta1} and \eqref{def:mu}, we obtain that
$$\begin{aligned}
(II_1)^{\frac{1}{p_1}} &\le c r^{\alpha}\left(\fiint_{Q_r}(1+|D\bfu|)^{p}\,\mathrm{d}z\right)^{\frac{q-1}{p}}\\
&\le c r^{\alpha}\left(\fiint_{Q_r}H^-(1+|D\bfu|)\,\mathrm{d}z\right)^{\frac{q-1}{p}}\\
&\le c r^{\alpha}\left(\fiint_{Q_r}H^-_{1+|D\tilde{\bm\ell}|}(|D\bfu_{\tilde{\bm\ell}}|)\,\mathrm{d}z + H^-(1+|D\tilde{\bm\ell}|)\right)^{\frac{q-1}{p}} \le c r^{\alpha} \le c \mu.
\end{aligned}$$
Similarly, in addition, using the Jensen type inequality in Lemma~\ref{lem:Jensen} with $\psi(s)=H^{-,*}_{\sigma}(s^{1/p_1})$ and using $H^{-,*}_{\sigma}\circ (H^-_{\sigma})'\sim H^-_{\sigma}$, we also obtain that
$$\begin{aligned}
(II_2)^{\frac{1}{p_1}} & \le c r^{\beta_0} (H^{-,*})^{-1}\left(\fiint_{Q_r} H^-(1+|D\bfu|)\,\d z\right)\\
&\le c r^{\beta_0} (H^{-,*})^{-1}\left(\fiint_{Q_r}H^-_{1+|D\tilde{\bm\ell}|}(|D\bfu_{\tilde{\bm\ell}}|)\,\mathrm{d}z + H^-(1+|D\tilde{\bm\ell}|)\right) \le c r^{\beta_0} \le c \mu\,,
\end{aligned}$$
and 
$$\begin{aligned}
(II_3)^{\frac{1}{p_1}} \le (H^{-,*}_{1+|D\tilde{\bm\ell} |})^{-1} \left(\fiint_{Q_r} H^-_{1+|D\tilde{\bm\ell} |}(|D \bfu_{\tilde{\bm\ell}}|)\, \d z\right)\le c (H^{-,*}_{1+|D\tilde{\bm\ell} |})^{-1} (\mu^2)\,.
\end{aligned}$$
Therefore,
$$
\left(\fiint_{Q_r} |{\bf G}|^{p_1}\, \d z \right)^{\frac{1}{p_1}}\le c \left(\mu+ (H^{-,*}_{1+|D\tilde{\bm\ell} |})^{-1} (\mu^2)\right)\,.
$$
Moreover, using the Jensen type inequality in Lemma~\ref{lem:Jensen} with $\psi(s)=H^{-}_{\sigma}(s^{1/p_1})$, we also obtain 
$$
\left(\fiint_{Q_r} |D \bfu_{\tilde{\bm\ell}}|^{p_1}\, \d z\right)^{\frac{1}{p_1}}  \le c (H^-_{1+|D\tilde{\bm\ell}|})^{-1} \left(\fiint_{Q_r} H^-_{1+|D\tilde{\bm\ell}|} (|D \bfu_{\tilde{\bm\ell}}|)\, \d z \right) \le c (H^-_{1+|D\tilde{\bm\ell}|})^{-1} (\mu^2)\,.
$$
We further observe from \eqref{Dlequiv}, \eqref{(2.6a)} and \eqref{(2.6b)} that for any  $s\in(0 ,1]$,
$$
(H^\pm_{1+|D\tilde{\bm\ell}|})'(s) = \frac{(H^\pm)'(1+|D\tilde{\bm\ell}|+s)}{1+|D\tilde{\bm\ell}|+s} s \sim s\,,
$$
which together with \eqref{(2.6a)} and \eqref{eq:hok2.41} implies 
\begin{equation}\label{sequiv}
s^2 \sim  H^\pm_{1+|D \tilde{\bm\ell}|} (s)
\ \ \text{and}\ \
s^2 \sim ((H^\pm_{1+|D\tilde{\bm\ell}|})')^{-1}(s) s \sim (H^{\pm,*}_{1+|D\tilde{\bm\ell}|})'(s) s \sim  H^{\pm,*}_{1+|D\tilde{\bm\ell}|}(s) \,.
\end{equation}
Finally, applying \eqref{sequiv} to $s=\mu$, we obtain
\begin{equation}\label{Dup1}
\left(\fiint_{Q_r} |D \bfu_{\tilde{\bm\ell}}|^{p_1}\, \d z+ \fiint_{Q_r} |{\bf G}|^{p_1}\, \d z\right)^{\frac{1}{p_1}}
\le c \mu \,.
\end{equation}

Therefore, since the estimate \eqref{Acalori_pf1} with $\bm\ell=\tilde{\bm\ell}$, together with  \eqref{Dlequiv}, \eqref{def:mu} and \eqref{sequiv} for $s=\mu$, yields
$$
\begin{aligned}
&\frac{1}{|Q_r|}\left|\iint_{Q_r} \bfu_{\tilde{\bm\ell}}  \cdot \bm\zeta_t - \langle D_{\bm\xi}\bA(z^-,D\tilde{\bm\ell})  D\bfu_{\tilde{\bm\ell}} \,|\, D\bm\zeta\rangle \, \d z  - \left[\int_{B_r} \bfu_{\tilde{\bm\ell}}  \cdot \bm\zeta  \, \d x \right]^{t=r^2}_{t=-r^2}  \right| \\
& \le c\left( r^{\frac{\beta_1}{2}}+\delta_0^{\frac{\gamma}{2}}\right)  \mu \|D\bm\zeta\|_{L^\infty(B_{r},\R^{N\times n})}\,,
\end{aligned}
$$
we may apply Theorem~\ref{thm:Acaloric} to $\A:=D_{\bm\xi}\bA(z^-,D\tilde{\bm\ell})$ (this satisfies the Legendre-Hadamard condition by \eqref{eq:1.8ok1} and \eqref{eq:1.8ok2} and \eqref{Dlequiv}), $\psi(s):=s^{p_0}$ and $\sigma_0:=\tfrac{p_1}{p_0}-1$, where $p_1$ and $p_0$ are given in \eqref{def:p0p1}.
Thus,
for $\epsilon\in(0,1)$ to be determined small later, there exist small $\delta_0,r>0$ depending on $n,N,L,\nu,p,q,\gamma,K$ and $\epsilon$ such that  
\begin{equation}\label{comparisonp0}
\fiint_{Q_{r}} |D \bfu_{\tilde{\bm\ell}}- D \bfh|^{p_0}\, \d z \le \epsilon \mu^{p_0}\,,
\end{equation}
where $\bfh$ is the weak solution to 
\begin{equation}\label{eq:h}
\begin{cases}
\partial_t  \bfh- \div (\A D\bfh) = {\bf0}  \quad \text{in}\ \ Q_r, \\
 \bfh =\bfu_{\tilde{\bm\ell}} \quad \text{on}\ \ \partial_{\mathrm{p}} Q_r\,.
\end{cases}\end{equation}

Note that by the Lipschitz estimate \eqref{DhLip}, \eqref{comparisonp0}, \eqref{Dup1}, and \eqref{sequiv}
we have that
\begin{equation}\label{Dhphisigma}
\begin{aligned}
\fiint_{Q_{r/2}}H^-_{1+|D \tilde{\bm\ell}|} (|D\bfh|) \,\d z 
& \le c H^-_{1+|D \tilde{\bm\ell}|} \left(\fiint_{Q_{r}}|D{\bf h}| \,\d z \right)\\
& \le c H^-_{1+|D \tilde{\bm\ell}|} \left(\bigg(\fiint_{Q_{r}}|D{\bfu}_{\tilde{\bm\ell}}|^{p_0} \,\d z\bigg)^{\frac{1}{p_0}} +\mu\right)\\
 & \le  c H^-_{1+|D\tilde{\bm\ell}|} (\mu)  \le c \mu^2\,.
\end{aligned}\end{equation}
Therefore, recalling $1<\tilde p\le \tilde q$ given in \eqref{phiapq} and $\theta\in (0,1)$ given in Lemma~\ref{thm:sob-poincare} and additionally assumed that $\tilde q \theta>1$ and $\tilde p \theta > 2n/(n+2)$, letting $\kappa_0\in(0,1)$ such that $\theta = \tfrac{1-\kappa_0}{\tilde q} + \kappa_0$ (i.e., $\kappa_0=\tfrac{\tilde q \theta-1}{\tilde q-\theta}$),
and using H\"older's inequality, the Jensen type inequality in Lemma~\ref{lem:Jensen} with $\psi^{-1}(s):=[H^-_{1+|D \tilde{\bm\ell}|}(s)]^{1/\tilde q}$, the estimates 
 \eqref{comparisonp0} and \eqref{Dhphisigma}, \eqref{def:mu}, and \eqref{sequiv}, we obtain
\begin{equation}\label{comparisonphi}\begin{aligned}
&\left(\fiint_{Q_{r/2}} \left[H^-_{1+|D \tilde{\bm\ell}|}(|D \bfu_{\tilde{\bm\ell}}- D \bfh|)\right]^{\theta}\, \d z\right)^{\frac1\theta}\\ 
& \le \left(\fiint_{Q_{r/2}} \left[H^-_{1+|D\tilde{\bm\ell}|}(|D \bfu_{\tilde{\bm\ell}}- D \bfh |)\right]^{\frac{1}{\tilde q}}\, \d z\right)^{\frac{1-\kappa_0}{\theta}}\left(\fiint_{Q_{r/2}} H^-_{1+|D \tilde{\bm\ell}|}(|D\bfu_{\tilde{\bm\ell}}- D \bfh|)\, \d z\right)^{\frac{\kappa_0}{\theta}}\\
& \le c \left[H^-_{1+|D\tilde{\bm\ell}|}\left(\left(\fiint_{Q_{r/2}} |D \bfu_{\tilde{\bm\ell}}- D \bfh|^{p_0}\, \d z\right)^{\frac{1}{p_0}}\right)\right]^{\frac{1-\kappa_0}{\tilde{q}\theta}} \mu^{\frac{2\kappa_0}{\theta}}\\
& \le c \epsilon^{\frac{\tilde p(1-\kappa_0)}{p_0\tilde{q}\theta}} H^-_{1+|D\tilde{\bm\ell}|}(\mu)^{\frac{1-\kappa_0}{\tilde{q}\theta}} \mu^{2\frac{\kappa_0}{\theta}} 
\le c  \epsilon^{\frac{\tilde p (1-\kappa_0)}{p_0\tilde{q}\theta}}  \mu^{2}\,.
\end{aligned}\end{equation}
Therefore, by Lemma~\ref{thm:sob-poincare} and \eqref{def:mu} with \eqref{def:beta1},
 $$\begin{aligned}
\fiint_{Q_{r/2}} H^+_{1+|D \tilde{\bm\ell}|}\left(\frac{| \bfu_{ \tilde{\bm\ell}}-  \bfh|}{r}\right) \, \d z 
& \le c   \left(\fiint_{Q_{r/2}} \left[H^-_{1+|D \tilde{\bm\ell}|}(|D \bfu_{ \tilde{\bm\ell}}- D \bfh|)\right]^{\theta}\, \d z\right)^{1/\theta} \\
&\qquad + c \left(r^{\alpha-\frac{(q-p)(n+2)}{p}}+ r^\alpha (1+|D\tilde{\bm\ell}|)^{q-p}\right)(1+|D\tilde{\bm\ell}|)^p\\
& \le c  \epsilon^{\frac{\tilde p (1-\kappa_0)}{p_0\tilde{q}\theta}}  \mu^{2} + r^{\alpha-\frac{(q-p)(n+2)}{p}}\\
& \le c  \big(\epsilon^{\frac{\tilde p(1-\kappa_0)}{p_0\tilde{q}\theta}}+ r^{\frac12[\alpha-\frac{(q-p)(n+2)}{p}]} \big)  \mu^{2} \,.
 \end{aligned}$$
 Moreover, by Sobolev inequality for each time $t$, the Jensen type inequality in Lemma~\ref{lem:Jensen} with $\psi(s)=H^-_{1+|D\tilde{\bm\ell|}}(s^{\frac{n+2}{2n}})$, since $p>\tfrac{2n}{n+2}$, and \eqref{comparisonphi} with \eqref{sequiv}, we also have 
 $$\begin{aligned}
 \fiint_{Q_{r/2}} \left|\frac{ \bfu_{\tilde{\bm\ell}}-  \bfh }{r}\right|^2 \, \d z
 & \le  
 c \left( \fiint_{Q_{r/2}} \left|D\bfu_{\tilde{\bm\ell}}- D \bfh \right|^{\frac{2n}{n+2}}\, \d z\right)^{\frac{n+2}{n}}\\
&\le c\left[(H^{-}_{1+|D\tilde{\bm\ell}|})^{-1}\left( \left(\fiint_{Q_{r/2}} [H^-_{1+|D\tilde{\bm\ell}|} (|D\bfu_{\tilde{\bm\ell}}- D \bfh |)]^{\theta}\, \d z\right)^{\frac{1}{\theta}}\right)\right]^{2}\\
&\le c\left[(H^{-}_{1+|D\tilde{\bm\ell}|})^{-1}\left( \epsilon^{\frac{\tilde p(1-\kappa_0)}{p_0\tilde{q}\theta}}  \mu^{2} \right)\right]^{2}\\
&\le c \epsilon^{\frac{2(\tilde p-1)(1-\kappa_0)}{p_0\tilde{q}\theta}} \mu^2\,.
\end{aligned}
$$
Therefore, collecting the last two estimates we obtain
\begin{equation}\label{comparisonuh}
\fiint_{Q_{r/2}} \left[ \left|\frac{ \bfu_{\tilde{\bm\ell}}-  \bfh }{r}\right|^2 + H^+_{1+|D \tilde{\bm\ell}|}\left(\frac{| \bfu_{\tilde{\bm\ell}}-  \bfh|}{r}\right) \right]\, \d z \le  c \tilde\epsilon \mu^2\,,
\end{equation}
where
$$
\tilde \epsilon:=\epsilon^{\frac{\tilde p(1-\kappa_0)}{p_0\tilde{q}\theta}}+ r^{\frac12[\alpha-\frac{(q-p)(n+2)}{p}]}+ \epsilon^{\frac{2(\tilde p-1)(1-\kappa_0)}{p_0\tilde{q}\theta}}\,.
$$

\textit{Step 3 (Decay estimates). }  Let $\tau\in (0,\textstyle\frac{1}{4})$ be a small constant determined later. We shall show that 
\begin{equation}\label{decay1}
\fiint_{Q_{2\tau r}}  \left|\frac{ \bfu_{ \bm\ell_{2{\tau} r}}}{2{\tau} r}\right|^2 \, \d z +
 \fiint_{Q_{2{\tau} r}}  H^+_{1+|D \bm\ell_{2{\tau} r}|}\left(\frac{| \bfu_{ \bm\ell_{2{\tau} r}}|}{2{\tau} r}\right) \, \d z \le c_2 \left(\tau^{-(n+2+\tilde q)}\tilde{\epsilon}+  \tau^2\right) \mu^2    \,,
\end{equation}
for some constant $c_2\ge 1$, where $\bm\ell_{2{\tau} r}:=\bm\ell_{z_0,2{\tau} r,\bfu}$ is the linear function defined in \eqref{eq:linearmap}.

Let $\bfh$ be the $\mathcal A$-caloric map of \eqref{eq:h}, and define
$$
\tilde{\bm\ell}_{\rho,\bfh} := \tilde{\bm\ell}_{z_0,\rho,\bfh} (x):= (D\bfh)_{Q_\rho(z_0)} (x-x_0) + \bfh (z_0), \quad z_0=(x_0,t_0)\,.
$$
Observe that for $z\in Q_{r/2}$,
$$
\bfh(z)-\bfh(z_0)=\left(\int_0^1D\bfh(s z + (1-s)z_0) \, \mathrm{d}s \right) (x-x_0) +   \left(\int_0^1\partial_t\bfh(s z + (1-s)z_0) \, \mathrm{d}s \right) (t-t_0)\,,
$$
from which together with \eqref{supestimate} with $\bfh$ in place of $D_i\bfh$, $i=1,2,\dots,n$, \eqref{DhLip}, \eqref{Dup1} and \eqref{comparisonp0} 
\begin{equation}\label{hzhz0}
\begin{aligned}
\bfh(z)-\bfh(z_0)
& \le \frac{r}{2}\|D\bfh\|_{L^\infty(Q_{r/2})} +  \frac{r^2}{4}\|\bfh_t\|_{L^\infty(Q_{r/2})}\\
& \le c \left(r\|D\bfh\|_{L^\infty(Q_{r/2})} +  r^2\|D^2\bfh \|_{L^\infty(Q_{r/2})}\right)\\
&\le c r \fint_{Q_{r}}|D\bfh|\, \d z \le c r \left(\fint_{Q_{r}}|D\bfh-D\bfu_{\tilde{\bm\ell}}|\, \d z+\fint_{Q_{r}}|D\bfu_{\tilde{\bm\ell}}|\, \d z\right)  \le c r \mu.
\end{aligned}\end{equation}
Similarly, we also obtain that 
\begin{equation}
\begin{aligned}
\sup_{Q_{2{\tau} r}} |\bfh - \tilde{\bm\ell}_{2{\tau} r, \bfh}| &\le   2{\tau} r  \sup_{Q_{2{\tau} r}} |D\bfh - (D\bfh )_{Q_{2{\tau} r}}|+ (2{\tau} r)^2 \sup_{Q_{2{\tau} r}} |\partial_t \bfh|\\
& \le  (2{\tau} r)^2 \sup_{Q_{2{\tau} r}}|D^2\bfh|  +(2{\tau} r)^3 \sup_{Q_{2{\tau} r}}|D\partial_t \bfh| + (2{\tau} r)^2 \sup_{Q_{2{\tau} r}}|\partial_t \bfh|\\
& \le c\tau^2 r  \left( r  \sup_{Q_{r/2}} |D^2\bfh|  + r^2 \sup_{Q_{r/2}} |D^3 \bfh|\right)\\
& \le c \tau^2 r \left(\fiint_{Q_{r}} |D\bfh|^{p_0}\,\mathrm{d}z \right)^{1/{p_0}}\le  c \tau^2 r \mu\,.
\end{aligned}
\label{eq:estimateh2}
\end{equation}

Let us first consider the second integral in \eqref{decay1}.  Observe  that 
by  \eqref{DlDl} with $\rho=2{\tau} r$ and $\bm\ell= \tilde{\bm\ell} +{\bf h}(z_0)$, $\eqref{comparisonuh}$, \eqref{hzhz0} and \eqref{def:mu},
 $$\begin{aligned}
|D\bm\ell_{2{\tau} r}| 
&\le |D\tilde{\bm\ell}| + (n+2) \fiint_{Q_{2{\tau} r}} \frac{|\bfu_{\tilde{\bm\ell}}-{\bf h}(z_0)|}{2{\tau} r}\,\mathrm{d}z \\
&\le |D\tilde{\bm\ell}| + \frac{c}{\tau^{n+3}} \left(\fiint_{Q_{r/2}} \left|\frac{\bfu_{\tilde{\bm\ell}} -{\bf h}}{r}\right|^2\,\mathrm{d}z+\fiint_{Q_{r/2}} \left|\frac{{\bf h} -{\bf h}(z_0)}{r}\right|^2\,\mathrm{d}z\right)^{1/2} \\
&\le |D\tilde{\bm\ell}| +  c \tau^{-n-3} \mu \\
&\le |D\tilde{\bm\ell}| +  c \tau^{-n-3} \delta_0^{\frac{1}{2}} \le K+1 \,. 
\end{aligned}$$
In the last inequality, we have chosen  $\delta_0>0$ sufficiently small such that $c(\tfrac{\tau}{2})^{-n-3}\delta_0^{1/2}\le 1$ with $\tau>0$ that will be determined in the end of the proof.  
Then it follows that $1+|D\bm\ell_{2{\tau} r}|\sim 1 \sim 1+ |D\tilde{\bm\ell}|$, and hence  $H^+_{1+|D\bm\ell_{2{\tau} r}|}(s)\sim H^+_{1+|D\tilde{\bm\ell}|}(s)$. Moreover,
by \eqref{eq:filv1} with $\varphi=H^+_{1+|D\tilde{\bm\ell}|}$ and $\bm\ell= \tilde{\bm\ell}+\tilde{\bm\ell}_{2{\tau} r, \bfh}$, 
\begin{equation}\begin{aligned}
& \fiint_{Q_{2{\tau} r}}  H^+_{1+|D \bm\ell_{2{\tau} r}|}\left(\frac{| \bfu_{ \bm\ell_{2{\tau} r}}|}{2{\tau} r}\right) \, \d z
\le  c \fiint_{Q_{2{\tau} r}}  H^+_{1+|D \tilde{\bm\ell}|}\left(\frac{| \bfu_{ \bm\ell_{2{\tau} r}}|}{2{\tau} r}\right) \, \d z\\
&\le  c \fiint_{Q_{2{\tau} r}}  H^+_{1+|D \tilde{\bm\ell}|}\left(\frac{| \bfu_{\tilde{\bm\ell}}-\tilde{\bm\ell}_{2{\tau} r, \bfh}|}{2{\tau} r}\right) \, \d z\\
& \le  c \fiint_{Q_{2{\tau} r}}  H^+_{1+|D \tilde{\bm\ell}|}\left(\frac{| \bfu_{\tilde{\bm\ell}}-\bfh|}{2{\tau} r}\right) \, \d z + c \fiint_{Q_{2{\tau} r}}  H^+_{1+|D \tilde{\bm\ell}|}\left(\frac{| \bfh - \tilde{\bm\ell}_{2{\tau} r, \bfh}|}{2{\tau} r}\right) \, \d z\\
& \le  c \tau^{-(n+2+\tilde q)}\fiint_{Q_{r/2}}  H^+_{1+|D \tilde{\bm\ell}|}\left(\frac{| \bfu_{\tilde{\bm\ell}}-\bfh|}{ r}\right) \, \d z + c \fiint_{Q_{2{\tau} r}}  H^+_{1+|D\tilde{\bm\ell}|}\left(\frac{| \bfh - \tilde{\bm\ell}_{2{\tau} r, \bfh}|}{2{\tau} r}\right) \, \d z\,.
\end{aligned} 
\label{eq:estimateh1}
\end{equation}
Combining the estimates \eqref{eq:estimateh1}, \eqref{comparisonuh}, \eqref{eq:estimateh2} and \eqref{sequiv} yields
$$
\fiint_{Q_{2{\tau} r}}  H^+_{1+|D \bm\ell_{2{\tau} r}|}\left(\frac{| \bfu_{ \bm\ell_{2{\tau} r}}|}{2{\tau} r}\right) \, \d z 
\le   c \tau^{-(n+2+\tilde q)}\tilde{\epsilon}\mu^2 + c  H^+_{1+|D \tilde{\bm\ell}|}(\tau \mu)
\le c \tau^{-(n+2+\tilde q)}\tilde{\epsilon}\mu^2 + c \tau^2 \mu^2\,. 
$$
In a similar same way,  we can estimate the first integral in \eqref{decay1} as follows:
$$\begin{aligned}
\fiint_{Q_{2{\tau} r}}  \left|\frac{ \bfu_{ \bm\ell_{2{\tau} r}}}{2{\tau} r}\right|^2 \, \d z & \le c \fiint_{Q_{2{\tau} r}}  \left|\frac{ \bfu_{\tilde{\bm\ell}}-{\tilde{\bm\ell}_{2{\tau} r, \bfh}}}{2{\tau} r}\right|^2 \, \d z \\
& \le c  \fiint_{Q_{2{\tau} r}}  \left|\frac{ \bfu_{ \tilde{\bm\ell}}-\bfh}{2{\tau} r}\right|^2 \, \d z  + \fiint_{Q_{2{\tau} r}}  \left|\frac{ \bfh  -{\tilde{\bm\ell}_{2{\tau} r, \bfh}}}{2{\tau} r}\right|^2 \, \d z \\
& \le c \tau^{-(n+4)}\tilde{\epsilon}\mu^2 + c \tau^2\mu^2. 
\end{aligned}$$
Therefore, we obtain \eqref{decay1}. 

Now, we prove the inequality \eqref{nondegenerate_decay}. By \eqref{eq:excess}, \eqref{eq:filv3} for ${\bf Q}=D\bm\ell_{2\tau r}$, the Caccioppoli inequality \eqref{eq:caccio} for $\bm\ell = \bm\ell_{2\tau r}$, and \eqref{decay1}, we obtain
$$\begin{aligned}
\Phi(z_0, \tau r, \bfu) 
&=\fiint_{Q_{\tau r}} H^-_{1+|(D\bfu)_{z_0,\tau r}|}\left( \left|D\bfu -(D\bfu)_{z_0,\tau r}\right|\right) \, \d z + (\tau r)^{\frac{\beta_1}{2}} \\
& \le c\, \fiint_{Q_{\tau r}} H_{1+|D\bm\ell_{2\tau r}|}(z,|D\bfu_{\bm\ell_{2\tau r}}|)\, \d z +  (\tau r)^{\frac{\beta_1}{2}} \\
&\le c \left\{\fiint_{Q_{2\tau r}}  \left|\frac{ \bfu_{ \bm\ell_{2 \tau r}}}{2\tau r}\right|^2 + H^+_{1+|D \bm\ell_{2\tau r}|}\left(\frac{| \bfu_{ \bm\ell_{2\tau r}}|}{2\tau r}\right) \, \d z + (\tau r)^{\beta_1} + (\tau r)^{\frac{\beta1}{2}} \right\}  \\
&\le c_3 \left(\tau^{-(n+2+\tilde q)}\tilde{\epsilon}+  \tau^2\right) \mu^2 + c_3\tau^{\frac{\beta_1}{2}} r^{\frac{\beta_1}{2}}  \\
&\le \tau^{2\beta} \Phi(z_0,r, \bfu) \,,
\end{aligned}$$
where we chose $\tau\in(0, \tfrac14)$ and $\tilde \epsilon\in(0,1)$ such that 
$$
c_3 \tau^{\frac{\beta_1}{2}-2\beta}\le \frac{1}{2}\,, \quad 
4c_3 \tau^{2-2\beta} \le \frac{1}{4}\,,
\quad \text{and}\quad
 c_3 \tau^{-(n+2+\tilde q)} \tilde\epsilon \le \frac{1}{4}\tau^{2\beta}.
$$
(Then $\delta_0$ is also fixed.) This corresponds to \eqref{nondegenerate_decay}.

Finally, using  Jensen's inequality, \eqref{def:mu} and \eqref{sequiv}, we obatin 
$$\begin{aligned}
|(D\bfu)_{z_0,\tau r}| & \leq |(D\bfu)_{z_0,\tau r}-D\tilde{\bm\ell}| + |D\tilde{\bm\ell}| \\
& \leq \tau^{-(n+2)} \fiint_{Q_{r}} |D\bfu_{\tilde{\bm\ell}}|\,\mathrm{d}z + |D\tilde{\bm\ell}| \\
& \le c \tau^{-(n+2)} (H^-_{1+|D\tilde{\bm\ell}|})^{-1} \left(\fiint_{Q_r} H^-_{1+|D\tilde{\bm\ell}|} (|D\bfu_{\tilde{\bm\ell}}|)\, \d z \right) + |D\tilde{\bm\ell}| \\
& \le c \tau^{-(n+2)} (H^-_{1+|D\tilde{\bm\ell}|})^{-1} (\mu^2) + |D\tilde{\bm\ell}| \\
& \le c \tau^{-(n+2)} \Phi(z_0, r, \bfu)^{1/2} + |D\tilde{\bm\ell}| \\
& \le c_1 \tau^{-(n+2)} \delta_0^{1/2} + |D\tilde{\bm\ell}| \,
\end{aligned}
$$
for some constant $c_1\ge 0$, which implies \eqref{Dldecay}.
This completes the proof.
\end{proof}

The excess decay estimate of Lemma \ref{Lem:nondegenerate1} can be iterated on each scale, as expressed by the following lemma.
\begin{lemma}[Iteration] \label{Lem:iteration}
Let $\bfu$ be a weak solution to \eqref{system1}, $K\ge 1$, and $\beta\in(0,\tfrac{\beta_1}{4})$ where $\beta_1$ is given in \eqref{def:beta1}. There exist small $r_0\in (0,R_0]$ and $\delta_1,\tau_0 \in(0,1)$ depending on $n,N,p,q,\alpha,L,\nu,\gamma,K$ and $\beta$ such that if 
$Q_{r}(z_0)\Subset \Omega_T$ with $r\in(0,r_0]$,  
$$
1+|(D\bfu)_{z_0, r}| \le K \quad \mbox{ and } \quad \Phi(z_0,r,\bfu)  \le \delta_1\,, 
$$
then for each $j\in\mathbb{N}\cup\{0\}$ we have
\begin{itemize}
\item[(i)] $\Phi(z_0,\tau^jr,\bfu) \leq \tau^{2\beta j} \Phi(z_0,r,\bfu)\,,$ \\
\item[(ii)] $|(D\bfu)_{z_0,\tau^jr}| \leq K - \tau^{\beta j}\,,$
\end{itemize}
and moreover, for every $0<\rho\leq r$
\begin{equation}
\fiint_{Q_{\rho}(z_0)} H^-_{1+ |(D{\bfu})_{z_0,\rho}|}(|D\bfu - (D{\bfu})_{z_0,\rho} |)\,\mathrm{d}z \leq \bar{c} \left(\frac{\rho}{r}\right)^{2\beta} \Phi(z_0,r,\bfu)\,
\label{eq:filv7.8}
\end{equation}
for some $\bar{c}=\bar{c}(n,N,p,q,\alpha,L,\nu,\gamma,K,\beta)$.
\end{lemma}
\proof
Recall the constants $r_0$, $\delta_0$, $\tau$ and $c_1$ determined in  Lemma~\ref{Lem:nondegenerate1}, and set
\begin{equation}
\delta_1:= \left[\min \left\{ \delta_0,\, c_1^{-1}\tau^{n+2}  \big(1-\tau^{\beta}\big) \right\}\right]^{2} \,.
\label{eq:delta0small}
\end{equation}
Note that the inequalities in $(i)$ and $(ii)$ are trivial when $j=0$. Then, given $j_0\ge 0$, we assume that the inequalities in  $(i)$ and $(ii)$ hold for every $j=0,1,\dots,j_0$. By the definition of $\delta_1$ in \eqref{eq:delta0small}, $(i)$ and $(ii)$ for $j=j_0$, we deduce that
\begin{equation*}
|(D\bfu)_{z_0,\tau^{j_0}r}| \le K \quad \mbox{ and } \quad \Phi(z_0,\tau^{j_0}r,\bfu)\le \Phi(z_0,r,\bfu)\le \delta_1  \le {\delta}_0\,, 
\end{equation*}
so we can apply Lemma \ref{Lem:nondegenerate1} with $\tau^{j_0}r$ in place of $r$ to obtain
\begin{equation*}
\Phi(z_0,\tau^{j_0}r,\bfu) 
 \le  \tau^{2\beta} \Phi(z_0,\tau^{j_0+1}r,\bfu)
 \le \tau^{2\beta (j_0+1)} \Phi(z_0,r,\bfu)\,,
\end{equation*}
which proves $(i)$ for $j=j_0+1$. For what concerns $(ii)$, we have from  \eqref{Dldecay} with $\tau^{j_0}r$ in place of $r$ and $(i)$ for $j=j_0$ that
\begin{equation*}
\begin{split}
|(D\bfu)_{z_0,\tau^{j_0+1}r}| & \le  |(D\bfu)_{z_0,\tau^{j_0}r}| + c_1\tau^{-n-2} \Phi(z_0, \tau_{j_0}r, \bfu)^{\frac{1}{2}} \\
& \le K - \tau^{\beta j_0} + c_1\tau^{-n-2} \tau^{\beta j_0} \delta_1^{\frac12} \\
& \le K - \tau^{\beta j_0} + \tau^{\beta j_0} \left(1-\tau^{\beta}\right)  = K - \tau^{\beta(j_0+1)}\,,
\end{split}
\end{equation*}
which corresponds to $(ii)$ for $j=j_0+1$. Therefore, by induction, assertions $(i)$ and $(ii)$ hold for every integer $j\geq1$. We now prove \eqref{eq:filv7.8}. Let $\rho\in (0, r]$. There exists $j\geq0$ such that $\tau^{j+1}r < \rho \leq  \tau^{j}r$. Using \eqref{eq:filv3} 
and $(i)$ above, we can write
$$
\begin{aligned}
\fiint_{Q_{\rho}} H^-_{1+ |(D{\bfu})_{z_0,\rho}|}(|D\bfu - (D{\bfu})_{z_0,\rho} |)\,\mathrm{d}z  & \leq c \fiint_{Q_{\rho}} H^-_{1+  |(D{\bfu})_{z_0,\tau^jr}| }(|D\bfu - (D{\bfu})_{z_0,\tau^jr}|)\,\mathrm{d}z \\
& \leq \frac{c}{\tau^{n+2}} \fiint_{Q_{\tau^jr}} H^-_{1+ |(D{\bfu})_{z_0,\tau^jr}| }(|D\bfu - (D{\bfu})_{z_0,\tau_1^jr}|)\,\mathrm{d}z \\
& \leq  \frac{c}{\tau^{n+2}}\Phi(z_0,\tau^jr,\bfu) \\
& \leq  \frac{c}{\tau^{n+2}} \tau^{2\beta j} \Phi(z_0,r,\bfu) \\
& \leq \bar c \left(\frac{\rho}{r}\right)^{2\beta} \Phi(z_0,r,\bfu) \,,
\end{aligned}
$$
then the proof is concluded. 
\endproof

\section{Proof of Theorem \ref{thm:main-thm}} \label{sec:mainthm}

This section is entirely devoted to the proof of the main result, the partial regularity result of Theorem \ref{thm:main-thm}. 

\begin{proof}[Proof of Theorem \ref{thm:main-thm}]

Let $z_0\in\Omega_T$ be such that 
\begin{equation*}
\liminf_{r \to 0^+} \fiint_{Q_r(z_0)} |{\bf V}_{H^-_{Q_r(z_0)}}(D\bfu) - ({\bf V}_{H^-_{Q_r(z_0)}}(D\bfu))_{z_0,r}  |^2 \, \mathrm{d}z=0
\end{equation*} 
and 
\begin{equation*}
M_{z_0}:=\limsup_{r \to 0^+} | (D\bfu)_{z_0,r} |  < + \infty\,.
\end{equation*}
In particular, taking into account \eqref{monotonicity1H}, the first limit implies

$$
\mathop{\lim\inf}_{r\to0} \fiint_{Q_r} H^-_{1+|(D{\bfu})_{z_0,r}|}(|D{\bfu}-(D{\bfu})_{z_0,r}|) \,\mathrm{d}z =0\,, 
$$
where $H^-_{1+|(D{\bfu})_{z_0,r}|}$ is the shifted $N$-function of $H^-_{Q_r}$ with shift $1+|(D{\bfu})_{z_0,r}|$.
Let 
$$
K:= 2(2 + M_{z_0})\,,
$$ 
and, with $\delta_1$ and $r_0$ in Lemma \ref{Lem:iteration},
$$
\epsilon:= \frac{\delta_1}{4}\,, 
$$
and
$$
\bar{r}:= \min 
\left\{\left(\frac{\delta_1}{4}\right)^\frac{2}{\beta_1},\, r_0\right\} \,.
$$
Then we can find $r\leq \bar{r} $ small enough such that $Q_{2r}(z_0)\Subset \Omega_T$ and 
$$
\fiint_{Q_r(z_0)} H^-_{1+|(D{\bfu})_{z_0,r}|}(|D{\bfu}-(D{\bfu})_{z_0,r}|) \,\mathrm{d}z < \epsilon
\qquad \mbox{and} \qquad 
|(D\bfu)_{z_0,r} | < M_{z_0}+1\,,
$$
which implies that
$$
\Phi(z_0,r,\bfu) \leq \frac{\delta_1}{2} \quad \mbox{ and } \quad 1+|(D\bfu)_{z_0,r}| \leq \frac{K}{2} \,. 
$$
From the absolute continuity of the integrals, there exists $\tilde r\in (0, r]$ such that for every $\tilde z\in Q_{\tilde r}(z_0)\Subset \Omega_T$,
$$
\Phi(\tilde z,r,\bfu) \leq \delta_1 \quad \mbox{and} \quad 1+|(D\bfu)_{\tilde z,r}| \leq K \,. 
$$
Therefore, by applying Lemma \ref{Lem:iteration} to each $Q_r(\tilde z)\Subset\Omega_T$, in place of $Q_r(z_0)$,  with $\tilde z\in Q_{\tilde r}(z_0)$,  and recalling \eqref{monotonicity1H}, we deduce that for each $\beta\in(0,\tfrac{\beta_1}{4})$,
$$
\begin{aligned}
\fiint_{Q_{\rho}(\tilde z)} \frac{|\bV_p(D\bfu) - (\bV_p(D\bfu))_{\tilde z,\rho} |^2}{\rho^{2\beta}} \, \mathrm{d}z  &\le \frac{c}{\rho^{2\beta}} \fiint_{Q_{\rho}(\tilde z)} |{\bf V}_{H^-_{Q_{\rho}(\tilde z)}}(D\bfu) - ({\bf V}_{H^-_{Q_{\rho}(\tilde z)}}(D\bfu))_{\tilde z,\rho}  |^2 \, \mathrm{d}z \\
& \le \frac{c}{\rho^{2\beta}} \fiint_{Q_{\rho}(\tilde z)} H^-_{1+|(D{\bfu})_{\tilde z,\rho}|} (|D\bfu - (D{\bfu})_{\tilde z,\rho} |) \, \d z \\
& \leq  \frac{c\delta_1}{r^{2\beta}}  \,,
\end{aligned}
$$
for every $\rho \in (0, r]$ and $\tilde z\in Q_{\tilde r}(z_0)$. 
This implies $\bV_p(D\bfu) \in C^{0, \beta, \frac{\beta}{2}}(Q_{\tilde r}(z_0),\mathbb{R}^{N\times n})$, and concludes the proof. 
\end{proof}

\section*{Acknowledgments.}

J. Ok was supported by the National Research Foundation of Korea (NRF) funded by the Korean government (MSIT) (2022R1C1C1004523). 
The work of G. Scilla is part of the project “Variational Analysis of Complex Systems in Materials Science, Physics and Biology” PRIN Project 2022HKBF5C. J. Ok and B. Stroffolini were funded by the University of Naples Federico II through  the International Co-operation between University of Naples Federico II and Sogang University. The research of B. Stroffolini is part of the project
 “Geometric Evolution Problems and Shape
Optimizations”, PRIN  Project 2022E9CF89. The authors G. Scilla and B. Stroffolini are members of Gruppo Nazionale per l’Analisi Matematica, la
Probabilit\`a e le loro Applicazioni.
\bibliographystyle{amsplain}

\end{document}